\documentclass[11pt]{amsart}
\usepackage[latin1]{inputenc}
\usepackage[T1]{fontenc}
\usepackage{amsfonts}
\usepackage{amsmath}
\usepackage{amssymb}
\usepackage{graphicx}
\usepackage{pstricks}
\usepackage{pst-grad}
\usepackage{multido}
\usepackage{amsthm}
\usepackage{geometry}
\newcommand{\cA}{\mathcal{A}}
\newcommand{\cB}{\mathcal{B}}
\newcommand{\cC}{\mathcal{C}}

\newcommand{\cH}{\mathcal{H}}

\newcommand{\cL}{\mathcal{L}}
\newcommand{\cM}{\mathcal{M}}

\newcommand{\cR}{\mathcal{R}}

\newcommand{\cU}{\mathcal{U}}

\newcommand{\sg}{\mathcal{h}}
\newcommand{\sd}{\mathcal{i}}

\newcommand{\fA}{\mathfrak{A}}
\newcommand{\fB}{\mathfrak{B}}

\newcommand{\nZ}{\mathbb{Z}}

\newcommand{\nN}{\mathbb{N}}

\newcommand{\tC}{\tilde{C}}

\newcommand{\tG}{\tilde{G}}

\newcommand{\sqs}{\sqsubset}
\newcommand{\sq}{\sqsubseteq}
\newcommand{\ov}{\overline}
\newcommand{\ind}{\underset}

\newcommand{\tx}{\quad\text}

\newtheorem{Th}{Theorem}[section]

\newtheorem{Lem}[Th]{Lemma}
\newtheorem{Cor}[Th]{Corollary}
\newtheorem{Prop}[Th]{Proposition}
\newtheorem{Def-Prop}[Th]{Definition-Proposition}

\theoremstyle{definition}
\newtheorem{Def}[Th]{Definition}
\newtheorem{Exa}[Th]{Example}
\newtheorem{Cl}[Th]{Claim}
\theoremstyle{remark}
\newtheorem{Rem}[Th]{Remark}

\begin{document}

\bibliographystyle{plain}
\title{Generalized Induction of Kazhdan-Lusztig cells}
\author{J\'er\'emie Guilhot}
\address{Department of Mathematical Sciences, King's College, Aberdeen University,
Aberdeen AB24 3UE, Scotland, U.K.}
\address{Universit\'e de Lyon, Universit\'e Lyon 1, 
Institut Camille Jordan, CNRS UMR 5208, 
43 Boulevard du 11 Novembre 1918, F-69622 Villeurbanne Cedex, France.}
\curraddr{School of Mathematics and Statistics F07
University of Sydney NSW 2006
Australia }
\email{guilhot@maths.usyd.edu.au}

\date{November, 2007}
\begin{abstract}
Following Lusztig, we consider a Coxeter group $W$ together with a weight function. Geck showed that the Kazhdan-Lusztig cells of $W$ are compatible with parabolic subgroups. In this paper, we generalize this argument to some subsets of $W$ which may not be parabolic subgroups. We obtain two applications: we show that under specific technical conditions on the parameters, the cells of a certain parabolic subgroup of $W$ are cells in the whole group, and we decompose the affine Weyl group $\tilde{G}_{2}$ into left and two-sided cells for a whole class of weight functions.

\end{abstract}

\maketitle


\section{Introduction}
This paper is concerned with the partition of  a Coxeter group $W$ (more specifically affine Weyl groups) into Kazhdan-Lusztig cells with respect to a weight function, following the general setting of Lusztig \cite{bible}. This is known to play an important role in the representation theory of the corresponding Hecke algebra, Lie algebra and group of Lie type.

In the case where $W$ is an integral and bounded Coxeter group (see \cite[Chap. 1]{bible}) and $L$ is constant on the generators of $W$ (equal parameter case), there is an interpretation of the Kazhdan-Lusztig polynomials in terms of intersection cohomology (see \cite{KL}) which leeds to many deep properties, for which no elementary proofs are known. For instance, the coefficients of the Kazhdan-Lusztig polynomials are non-negative integers. In that case, the left cells have been explicitly described for the affine Weyl groups of type $\tilde{A}_{r},r\in\mathbb{N}$ (see \cite{Lus2, Shi1}), ranks 2, 3 (see \cite{Bed, Du, Lus3}) and types $\tilde{B}_{4}$, $\tilde{C}_{4}$ and $\tilde{D}_{4}$ (see \cite{Chen,Shi3,Shi2}).

Much less is known for unequal parameters. Lusztig has formulated a number of precise conjectures in that case (see \cite[\S 14, P1-P15]{bible}). The left cells have been explicitly described for the affine Weyl groups of type $\tilde{A}_{1}$ for any parameters (\cite{bible}) and $\tilde{B}_{2}$ when the parameters are coming from a graph automorphism (\cite{Bremke}). Note that the proof in the case $\tilde{B}_{2}$ involved the positivity property of the Kazhdan-Lusztig polynomials in the equal parameter case.

One of the few things which are known in the general case of unequal parameters, is the compatibility of the left cells with parabolic subgroups; see \cite{geck}. In a precise sense, any left cell of a parabolic subgroup can be ``induced'' to obtain a union of left cells of the whole group $W$. The main observation of this paper is that the methods of \cite{geck} work in somewhat more general settings, so that we can ``induce'' from subsets of $W$ which are not parabolic subgroups (see Section 3). This leads to our ``Generalized Induction Theorem''. 

We discuss two applications of this theorem. First we show the following result; see Section~4.
\begin{Th}
\label{main4}
Let $(W,S)$ be an arbitrary Coxeter system together with a weight function $L$. Let $W'\subseteq W$ be a  bounded standard parabolic subgroup with generating set $S'$ and let $N\in\nN$ be a bound for $W'$.  If $L(t)>N$ for all $t\in S-S'$ then the left cells (resp. two-sided cells) of $W'$, considered as a proper Coxeter group, are left cells (resp. two-sided cells) of $W$. 
\end{Th}

Then, we decompose the affine Weyl groups $\tG_{2}$ into left and two-sided cells for a whole class of weight functions. Namely, the ones which satisfy $L(s_{1})>4L(s_{2})=4L(s_{3})$ where 
$$\tG_{2}:=\sg s_{1},s_{2},s_{3}\ |\ (s_{1}s_{2})^{6}=1,(s_{2}s_{3})^{3}=1,(s_{1}s_{3})^{2}=1\sd.$$
We also determine the partial left (resp. two-sided) order on left (resp. two-sided) cells; see Section 6.


\section{Hecke algebra and geometric realization of an affine Weyl group}
\subsection{Hecke algebra and Kazhdan-Lusztig cells}
In this section, $(W,S)$ denotes an arbitrary Coxeter system. The basic reference is \cite{bible}. Let $L$ be a weight function. Recall that a weight function on $W$ is a function $L:W\rightarrow \mathbb{Z}$ such that $L(ww')=L(w)+L(w')$ whenever $\ell(ww')=\ell(w)+\ell(w')$. In this paper, we shall only consider the case where $L(w)>0$ for all $w\neq e$ (where $e$ is the identity element of $W$). A weight function is completely determined by its values on $S$ and must only satisfy $L(s)=L(t)$ if $s,t\in S$ are conjugate.

Let $\cA=\nZ[v,v^{-1}]$ and $\cH$ be the Iwahori-Hecke algebra corresponding to $(W,S)$ with parameters $\{L(s)\mid s\in S\}$. Thus $\cH$ has an $\cA$-basis $\{T_{w}\mid w\in W\}$, called the standard basis, with multiplication given by
\begin{equation*}
T_{s}T_{w}=
\begin{cases}
T_{sw}, & \mbox{if } sw>w,\\
T_{sw}+(v^{L(s)}-v^{-L(s)})T_{w}, &\mbox{if } sw<w,
\end{cases}
\end{equation*}
(here, ``<'' denotes the Bruhat order) where $s\in S$ and $w\in W$. \\
Let $\cA_{<0}=v^{-1}\nZ[v^{-1}]$ and $\cA_{\leq 0}=\nZ[v^{-1}]$. For $x,y\in W$ we set
$$T_{x}T_{y}=\sum_{z\in W} f_{x,y,z}T_{z}\quad\text{where $f_{x,y,z}\in \cA$}.$$
We say that $N\in\nN$ is a bound for $W$ if $v^{-N}f_{x,y,z}\in\cA_{\leq 0}$ for all $x,y,z$ in $W$. If there exists such a $N$, we say that $W$ is bounded. \\
Let $a\mapsto \overline{a}$ be the involution of $\cA$ which takes $v^{n}$ to $v^{-n}$ for all $n\in\nZ$.
We can extend it to a ring involution from $\cH$ to itself with
$$\overline{\underset{w\in W}{\sum}a_{w}T_{w}}=\underset{w\in W}{\sum}\overline{a}_{w}T^{-1}_{w^{-1}}\ , \text{ where $a_{w}\in \cA$}.$$

For $w\in W$ there exists a unique element $C_{w}\in\cH$ such that
$$\overline{C}_{w}=C_{w} \text{ and } C_{w}=T_{w}+\underset{y<w}{\underset{y\in W}{\sum}}P_{y,w}T_{w} $$
where $P_{y,w}\in \cA_{<0}$ for $y<w$. In fact, the set $\{C_{w},w\in W\}$ forms a basis of $\cH$, known as the Kazhdan-Lusztig basis. The elements $P_{y,w}$ are called the Kazhdan-Lusztig polynomials. We set $P_{w,w}=1$ for any $w\in W$.

Let $w\in W$ and $s\in S$, we have the following multiplication formula
$$C_{s}C_{w}=
\begin{cases}
C_{sw}+\underset{z;sz<z<w}{\sum}M_{z,w}^{s}C_{z}, & \mbox{if }w<sw,\\ 
(v_{s}+v_{s}^{-1})C_{w}, & \mbox{if }sw<w,
\end{cases}
$$
where $M_{z,w}^{s}\in \cA$ satisfies
\begin{equation*}
 \overline{M_{y,w}^{s}}=M_{y,w}^{s},
 \end{equation*}
 \begin{equation*}
 \label{M dependence}
(\underset{z;y\leq z<w;sz<z}{\sum}P_{y,z}M_{z,w}^{s})-v_{s}P_{y,w} \in \cA_{<0}.
\end{equation*}
It is shown in \cite[Proposition 6.4]{bible} that $M^{s}_{y,w}$ is a $\nZ$-linear combination of $v^{n}$ such that $-L(s)+1\leq n\leq L(s)-1$.\\ 
We have a similar formula for the multiplication on the right by $C_{s}$, we obtain some polynomials $M_{z,w}^{s,r}$ which satisfy $M_{z,w}^{s,r}=M_{z^{-1},w^{-1}}^{s}$.

The multiplication rule between the standard basis and the Kazhdan-Lusztig basis is as follows
\begin{equation*}
T_{s}C_{w}=
\begin{cases}
C_{sw}-v^{-L(s)}C_{w}+\underset{z;sz<z<w}{\sum}M_{z,w}^{s}C_{z}, & \mbox{if }w<sw,\\ 
v^{L(s)}C_{w}, & \mbox{if }sw<w.
\end{cases}
\end{equation*}

Let $y,w\in W$. We write $y\leftarrow_{L} w$ if there exists $s\in S$ such that $C_{y}$ appears with a non-zero coefficient in the expression of $T_{s}C_{w}$ (or equivalently $C_{s}C_{w}$) in the Kazhdan-Lusztig basis. The Kazhdan-Lusztig left pre-order $\leq_{L}$ on $W$ is the transitive closure of this relation. One can see that
$$\cH C_{w}\subseteq \underset{y\leq_{L} w}{\sum}\cA C_{y}\text{ for any $w\in W$.}$$
The equivalence relation associated to $\leq_{L}$ will be denoted by $\sim_{L}$ and the corresponding equivalence classes are called the left cells of $W$. Similarly, we define $\leq_{R}$, $\sim_{R}$ and right cells. We say that $x\leq_{LR} y$ if there exists a sequence
$$x=x_{0}, x_{1},..., x_{n}=y$$
such that for all $1\leq i\leq n$ we have $x_{i-1}\leftarrow_{L} x_{i}$ or $x_{i-1}\leftarrow_{R} x_{i}$. We write $\sim_{LR}$ for the associated equivalence relation and the equivalence classes are called two-sided cells. One can see that
$$\cH C_{w}\cH\subseteq \underset{y\leq_{LR} w}{\sum}\cA C_{y}\text{ for any $w\in W$.}$$
The pre-order $\leq_{L}$ (resp. $\leq_{LR}$) induces a partial order on the left (resp. two-sided) cells of $W$. 

For $w\in W$ we set $\cL(w)=\{s\in S| sw<w\}$ and $\cR(w)=\{s\in S| w>ws\}$.
It is shown in \cite[\S 8]{bible} that if $y\leq_{L} w$ then $\cR(w)\subseteq \cR(y)$. Similarly, if $y\leq_{R} w$ then $\cL(w)\subseteq \cL(y)$. 
We now introduce a definition.
\begin{Def}
Let $\fB$ be a subset of $W$. We say that $\fB$ is a left ideal of $W$ if and only if the $\cA$-submodule of $\cH$ generated by $\{C_{w}| w\in \fB\}$ is a left ideal. Similarly one can define right and two-sided ideals of $W$. 
\end{Def}
\begin{Rem}
 Here are some straightforward consequences of this definition\\
- Let $\fB$ be a left ideal and let $w\in \fB$. We have
$$\cH C_{w}\subseteq \sum_{y\in\fB}\cA C_{y}.$$
In particular, if $y\leq_{L} w$ then $y\in \fB$ and $\fB$ is a union of left cells.\\
- A union of left ideals of $W$ is a left ideal.\\
- An intersection of left ideals is a left ideal.\\
- A left ideal which is stable by taking the inverse is a two-sided ideal. In particular it is a union of two-sided cells.
\end{Rem}

\begin{Exa}
\label{Lid}
Let $J$ be a subset of $S$. We set
$$\cR^{J}:=\{w\in W\ |\ J\subseteq \cR(w)\}\quad\text{and}\quad \cL^{J}:=\{w\in W\ |\ J\subseteq \cL(w)\}.$$
Then the set $\cR^{J}$ is a left ideal of $W$. Indeed let $w\in \cR^{J}$ and $y\in W$ be such that $y\leq _{L} w$. Then we have $J\subseteq \cR(w)\subseteq \cR(y)$ and $y\in\cR^{J}$. Similarly one can see that $\cL^{J}:=\{w\in W| J\subseteq \cL(w)\}$ is a right ideal of $W$.
\end{Exa}

\subsection{A geometric realization}
\label{geom}
In this section, we present a geometric realization of an affine Weyl group. The basic references are \cite{Bremke, Lus1, Xi}. 

Let $V$ be an euclidean space of finite dimension $r\geq 1$. Let $\Phi$ be an irreducible root system of rank $r$ and $\check{\Phi}\subseteq V^{*}$ be the dual root system. We denote the coroot corresponding to $\alpha\in\Phi$ by $\check{\alpha}$ and we write $\mathcal{h}x,y\mathcal{i}$ for the value of $y\in V^{*}$ at $x\in V$. Fix a set of positive roots $\Phi^{+}\subseteq \Phi$. Let $W_{0}$ be the Weyl group of $\Phi$. For $\alpha\in\Phi^{+}$ and $n\in \mathbb{Z}$, we define a hyperplane
$$H_{\alpha,n}=\{x\in V\mid \mathcal{h}x,\check{\alpha}\mathcal{i}=n\}.$$
Let
$$\mathcal{F}=\{H_{\alpha,n}\mid \alpha\in \Phi^{+}, n\in\mathbb{Z}\}.$$
Any $H\in\mathcal{F}$ defines an orthogonal reflection $\sigma_{H}$ with fixed point set $H$. We denote by $\Omega$ the group generated by all these reflections, and we regard $\Omega$ as acting on the right on $V$. An alcove is a connected component of the set
$$V-\underset{H\in\mathcal{F}}{\bigcup}H.$$
$\Omega$ acts simply transitively on the set of alcoves $X$.

Let $S$ be the set of $\Omega$-orbits in the set of faces (codimension 1 facets) of alcoves. Then $S$ consists of $r+1$ elements which can be represented as the $r+1$ faces of an alcove. If a face $f$ is contained in the orbit $t\in S$, we say that $f$ is of type $t$.

Let $s\in S$. We define an involution $A\rightarrow sA$ of $X$ as follows.  Let $A\in X$; then $sA$ is the unique alcove distinct from $A$ which shares with $A$ a face of type $s$.  The set of such maps generates a group of permutations of $X$ which is a Coxeter group $(W,S)$. In our case, it is the affine Weyl group usually denoted $\tilde{W_{0}}$. We regard $W$ as acting on the left on $X$. It acts simply transitively and commutes with the action of $\Omega$.

Let $A_{0}$ be the fundamental alcove defined by
$$A_{0}=\{x\in V\mid 0<\sg x,\check{\alpha}\sd<1 \text{ for all $\alpha\in\Phi^{+}$}\}.$$

We illustrate this realization in Figure 1 in the case where $W$ is an affine Weyl group of type $\tilde{G}_{2}$ 
$$W:=\sg s_{1},s_{2},s_{3}\ |\ (s_{1}s_{2})^{6}=1,(s_{2}s_{3})^{3}=1,(s_{1}s_{3})^{2}=1\sd.$$
The thick arrows represent the set of positive roots $\Phi^{+}$, $zA_{0}$ and $yA_{0}$ are the image of the fundamental alcove $A_{0}$ under the action of $y=s_{2}s_{3}s_{2}s_{1}s_{2}s_{1}s_{2}\in W$ and $z=s_{3}s_{2}s_{1}s_{2}s_{1}s_{2}\in W$.

\begin{center}
\begin{pspicture}(-3,4.37)(3,-4.37)
\psset{unit=2.5cm}

\pslinewidth=.3pt


\pspolygon[fillcolor=lightgray,fillstyle=solid](0,0)(0,0.5)(0.2166,0.375)
\rput(0.085 ,0.35){{\footnotesize$A_{0}$}}
\psline[linewidth=0.5mm]{->}(0,0)(-1.3,0.75)
\psline[linewidth=0.5mm]{->}(0,0)(0.866,0)
\psline[linewidth=0.5mm]{->}(0,0)(-.433,0.75)
\psline[linewidth=0.5mm]{->}(0,0)(0,1.5)
\psline[linewidth=0.5mm]{->}(0,0)(0.433,0.75)
\psline[linewidth=0.5mm]{->}(0,0)(1.3,0.75)





\pspolygon[fillstyle=solid,fillcolor=lightgray](-0.433,-0.75)(-0.433,-0.25)(-0.2166,-0.375)
\rput(-0.33,-0.4){{\scriptsize $zA_{0}$}}

\pspolygon[fillstyle=solid,fillcolor=lightgray](-0.433,0.75)(-0.433,1.25)(-0.65,1.125)
\rput(-0.533,1.11){{\scriptsize $yA_{0}$}}

\psline(-1.732,1.5)(1.732,1.5)
\psline(-1.732,0.75)(1.732,0.75)
\psline(-1.732,0)(1.732,0)

\psline(-1.732,-0.75)(1.732,-0.75)
\psline(-1.732,-1.5)(1.732,-1.5)


\multido{\n=-1.732+0.433}{9}{
\psline(\n,-1.5)(\n,1.5)}

\psline(-1.732,1)(-0.866,1.5)
\psline(-1.732,0.5)(0,1.5)
\psline(-1.732,0)(0.866,1.5)
\psline(-1.732,-0.5)(1.732,1.5)

\psline(-1.732,-1)(1.732,1)

\psline(1.732,-1)(0.866,-1.5)
\psline(1.732,-0.5)(0,-1.5)
\psline(1.732,0)(-0.866,-1.5)
\psline(1.732,0.5)(-1.732,-1.5)


\psline(1.732,1)(0.866,1.5)
\psline(1.732,0.5)(0,1.5)
\psline(1.732,0)(-0.866,1.5)
\psline(1.732,-0.5)(-1.732,1.5)

\psline(1.732,-1)(-1.732,1)

\psline(-1.732,-1)(-0.866,-1.5)
\psline(-1.732,-0.5)(0,-1.5)
\psline(-1.732,0)(0.866,-1.5)
\psline(-1.732,0.5)(1.732,-1.5)

\psline(-1.732,0)(-0.866,1.5)
\psline(-1.732,-1.5)(0,1.5)
\psline(-0.866,-1.5)(0.866,1.5)
\psline(0,-1.5)(1.732,1.5)
\psline(0.866,-1.5)(1.732,0)

\psline(1.732,0)(0.866,1.5)
\psline(1.732,-1.5)(0,1.5)
\psline(0.866,-1.5)(-0.866,1.5)
\psline(0,-1.5)(-1.732,1.5)
\psline(-0.866,-1.5)(-1.732,0)
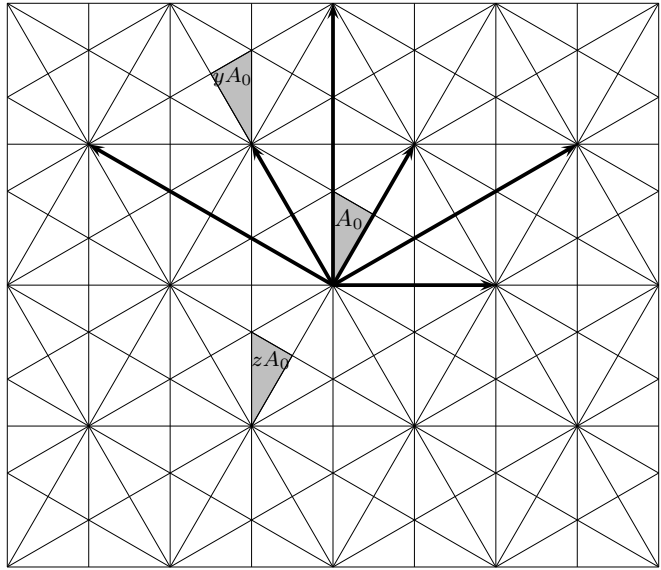
\rput(0,-1.75){Figure 1. Geometric realization of $\tilde{G}_{2}$}
\end{pspicture}
\end{center}


\section{Generalized induction of left cells}
\subsection{Main result}
\label{result}
Let $(W,S)$ be a Coxeter group together with a weight function $L$. Let $\cH$ be the associated Iwahori-Hecke algebra. In this section, we want to generalize the results of \cite{geck} on the induction of left cells. 

We consider a subset $U\subseteq W$ and a collection $\{X_{u}\ |\ u\in U\}$ of subsets of $W$ satisfying the following conditions
\begin{enumerate}
\item[{\bf I1}.] for all $u\in U$, we have $e\in X_{u}$,
\item[{\bf I2}.] for all $u\in U$ and $x\in X_{u}$ we have $\ell(xu)=\ell(x)+\ell(u)$,
\item[{\bf I3}.] for all $u,v\in U$ such that $u\neq v$ we have $X_{u}u\cap X_{v}v=\emptyset$,
\item[{\bf I4}.] the submodule $\cM:=\sg T_{x}C_{u}|\ u\in U,\ x\in X_{u}\sd_{\cA}\subseteq \cH$ is a left ideal,
\item[{\bf I5}.] for all $u\in U$, $x\in X_{u}$ and $u_{1}<u$ we have
$$P_{u_{1},u}T_{x}T_{u_{1}} \text{ is an $\cA_{<0}$-linear combination of $T_{z}$}.$$
\end{enumerate}

Let $u\in U$ and $x\in X_{u}$. We have
$$T_{x}C_{u}=T_{xu}+ \text{an $\cA$-linear combination of $T_{z}$ with $\ell(z)<\ell(xu)$}.$$
Since the set $\{T_{w}|w\in W\}$ is a basis of $\cH$, using {\bf I3}, one can see that $\cB=\{T_{x}C_{u}|u\in U,x\in X_{u}\}$ is a basis of $\cM$.\\

Let $u\in U$ and $z\in W$. Using {\bf I1}, {\bf I4} and the fact that $\cB$ is a basis of $\cM$, we can write
$$T_{z}C_{u}=\sum_{u\in U,x\in X_{u}}a_{x,u}T_{x}C_{u}\quad \text{for some $a_{x,u}\in \cA$}.$$
Let $\preceq$ be the relation on $U$ defined as follows. Let $u,v\in U$. We write $v\preceq u$ if there exist $x\in W$ and $z\in X_{v}$ such that $T_{z}C_{v}$ appears with a non-zero coefficient in the expression of $T_{x}C_{u}$ in the basis $\cB$. We still denote by $\preceq$ the pre-order induced by this relation (i.e the transitive closure). Since $C_{u}\in \cM$, we have
$$\cH C_{u}=\sum_{v\preceq u, z\in X_{v}}\cA T_{z}C_{v}.$$
\begin{Rem}
If we choose $U=W$ and $X_{w}=\{e\}$ for all $w\in W$, the pre-order $\preceq$ is the left pre-order $\leq_{L}$ on $W$.
\end{Rem}


We are now ready to state the main result of this section. 
\begin{Th}
\label{main}
Let $U$ be a subset of $W$ and $\{X_{u}|u\in U\}$ be a collection of subsets of $W$ satisfying conditions {\bf I1--I5}. Let $\mathcal{U}\subseteq U$ be such that
$$v\preceq u\in\cU \Longrightarrow v\in\cU.$$
Then, the set 
$$\{xu| u\in \cU,x\in X_{u}\}$$
is a left ideal of $W$.
\end{Th}
The proof of this theorem will be given in the next section. We have the following corollary.
\begin{Cor}
\label{mainC}
Let $\cC$ be an equivalence class on $U$ with respect to $\preceq$. Then the subset
$\{xu| u\in\cC, x\in X_{u}\}$ of $W$ is a union of left cells.
\end{Cor}
\begin{proof}
Let $v\in \cC$, $y\in X_{v}$ and $z\in W$ be such that $z\sim_{L} yv$. Consider the set $\cU=\{u\in U| u\preceq v\}$. Then $\cU$ satisfies the requirement of Theorem \ref{main}, thus $\fB:=\{xu| u\in \cU,\ x\in X_{u}\}$ is a left ideal of $W$. Since $z\leq_{L} yv$ and $yv\in \fB$, there exist $u_{z}\in \cU$ and $x\in X_{u_{z}}$ such that $z=xu_{z}$ and $u_{z}\preceq v$. \\
We also have $yv\leq_{L} xu_{z}$. Applying the same argument as above to the set $\{u\in U| u\preceq u_{z}\}$ yields that there exists $u_{y}\in \cU$ and $w\in X_{u_{y}}$ such that $yv=wu_{y}$ and $u_{y}\preceq u_{z}$. By condition {\bf I3}, we see that $u_{y}=v$. Thus $u_{z}\in \cC$ and the result follows.
\end{proof}
\begin{Rem} 
In \cite{geck}, Geck proved the following theorem, where $(W,S)$ is an arbitrary Coxeter system.
\begin{Th}
\label{parabolic}
Let $W'\subseteq W$ be a parabolic subgroup of $W$ and let $X'$ be the set of all $w\in W$ such that $w$ has minimal length in the coset $wW'$. Let $\cC$ be a left cell of $W'$. Then $X'\cC$ is a union of left cells of $W$.
\end{Th}
Let $U=W'$ and for all $w\in W'$ let $X_{w}=X'$. We claim that this theorem is a special case of Theorem \ref{main} and Corollary \ref{mainC}. Indeed, conditions {\bf I1--I3} and {\bf I5} are clearly satisfied. Condition {\bf I4} is a straightforward consequence of Deodhar's lemma; see \cite[Lemma 2.2]{geck}. Hence, it is sufficient to show that the pre-order $\preceq$ on $U=W'$ coincides with the Kazhdan-Lusztig left pre-order defined with respect to $W'$ (denoted $\leq'_{L}$) and the corresponding parabolic subalgebra $\cH':=\sg T_{w}\ |\ w\in W'\sd_{\cA}\subseteq \cH$. In other words, we need to show the following
$$u\leq'_{L} v \quad\Longleftrightarrow\quad u\preceq v.$$
Let $u,v\in W'$ such that $u\leq'_{L} v$. We may assume that there exists $s\in S'$ (where $S'$ is the generating set of $W'$) such that 
$$T_{s}C_{v}=\sum_{w\in W'} a_{w,v}C_{w} \quad \text{where $a_{w,v}\in\cA$ and $a_{u,v}\neq 0$}.$$
Since $C_{w}\in \cB$ for all $w\in W'$, this is the expression of $T_{s}C_{v}$ in $\cB$, which shows that $u\preceq v$.\\
Conversely, let $u,v\in W'$ such that $u\preceq v$. We may assume that there exist $z\in W$ and $x\in X'$ such that
$$T_{z}C_{v}=\sum_{w\in W',y\in X'}a_{yw,zv}T_{y}C_{w}\quad \text{where $a_{yw,zv}\in\cA$ and $a_{xu,zv}\neq 0$}.$$
We can write uniquely $z=z_{1}z_{0}$ where $z_{0}\in W'$, $z_{1}\in X'$ and $\ell(z)=\ell(z_{0})+\ell(z_{1})$. Then, we have
$$T_{z}C_{v}=T_{z_{1}}(T_{z_{0}}C_{v})=T_{z_{1}}(\sum_{w\in W', w\leq'_{L}v}a_{w,v}C_{w})=\sum_{w\in W', w\leq'_{L}v}a_{w,v}T_{z_{1}}C_{w}$$
and this is the expression of $T_{z}C_{v}$ in the basis $\cB$. We assumed that $T_{x}C_{u}$ appears with a non-zero coefficient, thus $u\leq'_{L} v$ as desired.
\end{Rem}
\subsection{Proof of Theorem \ref{main}}

We keep the setting of the last section and we introduce the following relation. Let $u,v\in U$, $x\in X_{u}$ and $y\in X_{v}$. We write $xu\sqs yv$ if $xu< yv$ (Bruhat order) and $u\preceq v$. We write $xu\sq yv$ if $xu\sqs yv$ or $x=y$ and $u=v$.

The main reference is the proof of \cite[Theorem 1]{geck}. 

\begin{Lem}
\label{rpol}
Let $v\in U$, $y\in X_{u}$. We have
$$T_{y^{-1}}^{-1}C_{v}=\sum_{u\in U,\ x\in X_{u}}\overline{r}_{xu,yv}T_{x}C_{u}$$
where $r_{yv,yv}=1$ and $r_{xu,yv}=0$ unless $xu\sq yv$.
\end{Lem}
\begin{proof}
Let $v\in U$ and $y\in X_{v}$. We have
$$T_{y^{-1}}^{-1}=T_{y}+\sum_{z<y}\ov{R}_{z,y}T_{z}$$
where $R_{z,y}\in\cA$ are the usual $R$-polynomials as defined in \cite[\S 4.3]{bible}. We obtain
\begin{align*}
T_{y^{-1}}^{-1}C_{v}&=(T_{y}+\sum_{z<y}\ov{R}_{z,y}T_{z})C_{v}\\
&=T_{y}C_{v}+\sum_{z<y}\ov{R}_{z,y}T_{z}C_{v}.
\end{align*}
Now we also have
$$T_{z}C_{v}=\ \text{$\cA$-linear combination of $T_{x}C_{u}$ where $u\preceq v$ and $x\in X_{u}$}.$$
We still have to show that if $T_{x}C_{u}$ appears in this sum then $xu<yv$.\\
This comes from the fact that $T_{z}C_{v}$, expressed in the standard basis, is an $\cA$-linear combination of term of the form $T_{w_{0}.w_{1}}$ where $w_{0}\leq z$ and $w_{1}\leq v$. In particular, since $z<y$ we have $w_{0}w_{1}<yv$. Then, expressing the right hand side of the equality in the standard basis, one can see that we must have $xu<yv$ if $T_{x}C_{u}$ appears with a non-zero coefficient.\\
Finally, by definition of $\sq$, we see that
$$T_{y^{-1}}^{-1}C_{v}=T_{y}C_{v}+\sum_{xu\sqs yv}\ov{r}_{xu,yv}T_{x}C_{u}.$$
The result follows.
\end{proof}
\begin{Lem}
\label{delta}
Let $u,v\in U$, $x\in X_{u}$ and $y\in X_{v}$. Then
$$\ind{xu\sq zw\sq yv}{\sum_{w\in U, z\in X_{w}}}\ov{r}_{xu,zw}r_{zw,yv}=\delta_{x,y}\delta_{u,v}$$

\end{Lem}
\begin{proof}
Since the map $h\mapsto \ov{h}$ is an involution and $C_{v}=\ov{C}_{v}$, we have
\begin{align*}
T_{y}C_{v}&=\ov{T_{y^{-1}}^{-1}C_{v}}\\
&=\ov{\sum_{w\in U,z\in X_{w}}\ov{r}_{zw,yv}T_{z}C_{w}}\\
&=\sum_{w\in U,z\in X_{w}}r_{zw,yv}T^{-1}_{z^{-1}}C_{w}\\
&=\sum_{w\in U,z\in X_{w}}r_{zw,yv}\big(\sum_{u\in U,x\in X_{u}}\ov{r}_{xu,zw}T_{x}C_{u}\big)\\
&=\sum_{u\in U,x\in X_{u}}\big(\sum_{w\in U,z\in X_{w}}\ov{r}_{xu,zw}r_{zw,yv}\big)T_{x}C_{u}.
\end{align*}
Since $\cB$ is a basis of $\cM$, using Lemma \ref{rpol} and comparing the coefficients yield the desired result.
\end{proof}

\begin{Prop}
\label{tilde}
Let $v\in U$ and $y\in X_{v}$. We have
$$C_{yv}=T_{y}C_{v}+\ind{xu\sqs yv}{\sum_{u\in U, x\in X_{u}}}p^{*}_{xu,yv}T_{x}C_{u}\qquad \text{where $p^{*}_{xu,yv}\in\cA_{<0}$}.$$
\end{Prop}
\begin{proof}
By Lemma \ref{delta}, there exists a unique family $(p^{*}_{xu,yv})_{xu\sqs yv}$ of polynomials in $\cA_{<0}$ such that 
$$\tC_{yv}:=T_{y}C_{v}+\ind{xu\sqs yv}{\sum_{u\in U, x\in X_{u}}}p^{*}_{xu,yv}T_{x}C_{u}$$
is stable under the $\ \bar{}\ $ involution; see \cite[p. 214]{duC0}, it contains a general setting to include the argument in \cite[Proposition 3.3]{geck} or in \cite[Theorem 5.2]{bible}.
Moreover, we have
\begin{align*}
\tilde{C}_{yv}&=T_{y}C_{v}+\ind{xu\sqs yv}{\sum_{u\in U, x\in X_{u}}}p^{*}_{xu,yv}T_{x}C_{u}\\
&=T_{y}\big( T_{v}+\sum_{v_{1}<v}P_{v_{1},v}T_{v_{1}} \big)+
\ind{xu\sqs yv}{\sum_{u\in U, x\in X_{u}}}p^{*}_{xu,yv}T_{x}\sum_{u_{1}\leq u}P_{u_{1},u}T_{u_{1}}\\ 
&=T_{yv}+\big(\sum_{v_{1}<v}P_{v_{1},v}T_{y}T_{v_{1}} \big)+
\ind{xu\sqs yv}{\sum_{u\in U, x\in X_{u}}}\sum_{u_{1}\leq u}p^{*}_{xu,yv}P_{u_{1},u}T_{x}T_{u_{1}}.\\
\end{align*}
By condition {\bf I5}, all the terms $P_{v_{1},v}T_{y}T_{v_{1}}$ occurring in the first sum and all the terms $p^{*}_{xu,yv}P_{u_{1},u}T_{x}T_{u_{1}}$ occurring in the second sum are $\cA_{<0}$-linear combinations of $T_{z}$ with $\ell(z)<\ell(yv)$. Thus
$$\tilde{C}_{yv}=T_{yv}+\ \text{an $\cA_{<0}$-linear combination of $T_{z}$ with $\ell(z)<\ell(yv)$}$$
and by definition and unicity of the Kazhdan-Lusztig basis, this implies that $\tilde{C}_{yv}=C_{yv}$.
\end{proof}

Let $\cU\subseteq U$ be as in Theorem \ref{main}. By definition of $\preceq$ one can see that 
$$\cM_{\cU}:=\sg T_{y}C_{v}\ |\ v\in \cU, \ y\in X_{v}\sd_{\cA}\subseteq \cH$$
is a left ideal. 
\begin{Cor}
\label{Cor1}
$$\cM_{\cU}=\sg C_{yv}\ |\ v\in \cU,\ y\in X_{v}\sd_{\cA}.$$
\end{Cor}
\begin{proof}
Let $v\in \cU$ and $y\in X_{v}$, using the previous proposition, we see that
$$C_{yv}=T_{y}C_{v}+\ind{xu\sqs yv}{\sum_{u\in \cU, x\in X_{u}}}p^{*}_{xu,yv}T_{x}C_{u}.$$
Thus $C_{yv}\in\cM_{\cU}$. Now, a straightforward induction on the order relation $\sq$ yields
$$T_{y}C_{v}=C_{yv}+\text{ an $\cA$-linear combination of $C_{xu}$}$$ 
where $u\in \cU,\ x\in X_{u}$ and $xu\sqs yv$.\\
This yields the desired assertion.
\end{proof}

We can now prove Theorem \ref{main}.\\
Let $\cU$ be a subset of $U$ such that
$$v\preceq u\in\cU \Longrightarrow v\in\cU.$$
Then $\cM_{\cU}=\sg T_{z}C_{w}\ |\ w\in \cU,\ z'\in X_{w}\sd_{\cA}\subseteq \cH$ is a left ideal. We want to show that the set $\fB:=\{yv|v\in \cU,\ y\in X_{v}\}$ is a left ideal of $W$.\\
Let $v\in \cU$, $y\in X_{v}$ and $z\in W$ be such that $z\leq_{L} yv$. We may assume that there exists $s\in S$ such that $C_{z}$ appears with a non-zero coefficient in the expression of $T_{s}C_{yv}$ in the Kazhdan-Lusztig basis. By Corollary \ref{Cor1}, we have $C_{yv}\in\cM_{\cU}$. Since $\cM_{\cU}$ is a left ideal we have $T_{s}C_{yv}\in\cM_{\cU}$. Thus, using Corollary \ref{Cor1} once more, we have
$$T_{s}C_{yv}=\sum_{u\in\cU,x\in X_{u}}a_{xu,yv}C_{xu}\quad \text{where $a_{xu,yv}\in \cA$}$$
and this is the expression of $T_{s}C_{yv}$ in the Kazhdan-Lusztig basis. The fact that $C_{z}$ appears with a non-zero coefficient in that expression implies that $z=xu$ for some $u\in\cU$ and $x\in X_{u}$. Thus $z\in \fB$, as desired. \hfill $\square$


\section{Cells in certain parabolic subgroups}
The aim of this section is to prove Theorem \ref{main4}. We will actually prove a stronger result.
Let $(W,S)$ be an arbitrary Coxeter system. For $J\subseteq S$, we denote by $X_{J}$ the set of  minimal left coset representatives with respect to the subgroup generated by $J$. Recall that $\cR^{J}=\{w\in W|\ J\subseteq \cR(w)\}$.
Let $W'\subseteq W$ be a standard parabolic subgroup with generating set $S'$. Furthermore, assume that $(W',S')$ is bounded by $N\in \nN$. 
\begin{Th}
\label{tfix}
Let $t\in S-S'$ be such that $L(t)>N$. Then 
$$\{w\in W|\ w=yw',\ y\in \cR^{\{t\}}\cap X_{S'}, w'\in W'\}$$
is a left ideal of $W$. 
\end{Th}
\begin{Rem}
This theorem implies Theorem \ref{main4}. Indeed, assume that, for all $t\in S-S'$ we have $L(t)>N$. Then
$$\bigcup_{t\in S-S'}\{w\in W|\ w=yw',\ y\in \cR^{\{t\}}\cap X_{S'}, w'\in W'\}=W-W'$$
is a left ideal of $W$. Furthermore, since it is stable by taking the inverse, it's a two-sided ideal. Thus $W-W'$ is a union of cells and so is $W'$. Let $y,w\in W'$ be such that $y\leq_{L} w$ in $W$. Then using Theorem \ref{parabolic}, one can easily see that $y\leq_{L} w$ in $W'$. Similarly, if $y\leq_{R} w$ in $W$ then $y\leq_{R} w$ in $W'$. The theorem follows.
\end{Rem}
Until the end of this section, we fix $t\in S-S'$ such that $L(t)>N$. \\
Let $U=tW'$. For $u\in U$ let 
$$X_{u}=(\cR^{\{t\}}\cap X_{S'})t .$$
We want to apply Theorem \ref{main} to the set $U$. One can directly check that conditions ${\bf I1}$--${\bf I3}$ hold. In order to check conditions ${\bf I4}$--${\bf I5}$ we need some preliminary lemmas. We denote by $\cH'$ the Hecke algebra associated to $(W',S')$ and the weight function $L$ (more precisely the restriction of $L$ to $S'$).

\begin{Lem}
\label{Cmult}
Let $w'\in W'$. We have
$$C_{t}C_{w'}=C_{tw'}\quad\text{ and }\quad T_{t}C_{w'}=C_{tw'}-v^{-L(t)}C_{w'}$$
\end{Lem}
\begin{proof}
We know that 
\begin{align*}
C_{t}C_{w'}&=C_{tw'}+\sum_{tz<z<w'}M_{z,w'}^{t}C_{z},\\
T_{t}C_{w'}&=C_{tw'}-v^{-L(t)}C_{w'}+\sum_{tz<z<w'}M_{z,w'}^{t}C_{z}.\\
\end{align*}
But $z<w'$ implies that $z\in W'$, thus we cannot have $tz<z$. The result follows.
\end{proof}
\begin{Rem}
Let $s'\in S'$. Since $L(t)\neq L(s')$, the order of $s't$ has to be even or infinite (otherwise, $s'$ and $t$ would be conjugate and $L(s')=L(t)$).
\end{Rem}
\begin{Lem}
\label{Tmult}
Let $s'\in S'$ and $w\in W'$. Let $m\in\nN$ be such that $m$ is less than or equal to the order of $s't$. We have
$$T_{(s't)^{m}}C_{w}=\sum_{w'\in W'}\sum^{m-1}_{i=0}a_{w',i}T_{(s't)^{i}s'}C_{tw'}+h'_{m}
$$
where $a_{w',i}\in\cA$ and $h'_{m}\in \cH'$, and
$$T_{(ts')^{m}}C_{w}=\sum_{w'\in W'}\sum^{m-1}_{i=0}b_{w',i}T_{(ts')^{i}}C_{tw'} +h''_{m}$$
where $b_{w',i}\in\cA$ and $h''_{m}\in \cH'$. Furthermore, $h'_{m}=h''_{m}$.
\end{Lem}
\begin{proof}
The first two equalities come from a straightforward induction.\\
It is clear that $h_{0}=h'_{0}=C_{w}$. Even though it is not necessary, let us do the case $m=1$ to show how the multiplication process works. We have
$$T_{s'}C_{w}=\sum_{w'\in W'} a_{w'}C_{w'}\quad \text{for some $a_{w'}\in\cA$.}$$
Thus we obtain (using the previous lemma)
\begin{align*}
T_{s't}C_{w'}&=T_{s'}C_{tw'}-v^{-L(t)}\sum_{w'\in W'} a_{w'}C_{w'}\\
\end{align*}
and
\begin{align*}
T_{ts'}C_{w'}&=\sum_{w'\in W'} a_{w'}C_{tw'}-v^{-L(t)}\sum_{w'\in W'} a_{w'}C_{w'}.\\
\end{align*}
It follows that 
$$h'_{1}=-v^{-L(t)}\sum_{w'\in W'} a_{w'}C_{w'}=h''_{1}.$$
Now, by induction, one can see that
$$h'_{m}=-v^{-L(t)}T_{s'}h'_{m-1}\in\cH'\quad\text{and}\quad h''_{m}=-v^{L(t)}T_{s'}h''_{m-1}\in\cH'.$$
The result follows.
\end{proof}
\begin{Prop}
\label{mod}
The submodule
$$\cM:=\sg T_{x}C_{u}\ |\ u\in U,\ x\in X_{u}\sd_{\cA}\subseteq \cH$$
is a left ideal.
\end{Prop}
\begin{proof}
Let $z\in W$, $u\in U$ and $x\in X_{u}$. We need to show that $T_{z}T_{x}C_{u}\in\cM$. Since $T_{z}T_{x}$ is an $\cA$-linear combination of $T_{y}$ ($y\in W$), it is enough to show that $T_{y}C_{u}\in \cM$ for all $y\in W$ and $u\in U$.\\
We proceed by induction on $\ell(y)$. If $\ell(y)=0$, then the result is clear.\\
Assume that $\ell(y)>0$. We may assume that $y\notin X_{u}$. Let $w'\in W'$ such that $u=tw'$. Recall that $X_{u}=(\cR^{\{t\}}\cap X_{S'})t$. \\
Suppose that $yt<y$, then we have
$$T_{y}C_{tw'}=T_{yt}T_{t}C_{tw'}=v^{L(t)}T_{yt}C_{tw'}\in \cM$$
by induction.\\
Suppose that $yt>y$. Since $yt\in \cR^{\{t\}}$ and $yt\notin \cR^{\{t\}}\cap X_{S'}$, there exists $s'\in S'$ such that $(yt)s'<yt$. Let $2n$ be the order of $ts'$ (it has to be finite in that case). One can see that there exists $y_{0}$ (with $\ell(y_{0})<\ell(y)$) such that $yt=y_{0}(ts')^{n}$. \\
Using Lemma \ref{Cmult} and the relation $C_{t}=T_{t}+v^{-L(t)}T_{e}$ we see that
$$C_{tw'}=C_{t}C_{w'}=T_{t}C_{w'}+v^{-L(t)}C_{w'}.$$
Since $s'\in S'$ and $w'\in W'$ we have
$$T_{s'}C_{w}=\sum_{w_{i}\in W'}a_{w_{i}}C_{w_{i}}\quad\text{for some $a_{w_{i}}\in \cA$}.$$
Thus we get
\begin{align*}
T_{y}C_{tw'}&=T_{yt}C_{w'}+v^{-L(t)}T_{y}C_{w'}\\
&=T_{y_{0}}T_{(ts')^{n}}C_{w'}+v^{-L(t)}T_{y_{0}}T_{(s't)^{n-1}s'}C_{w'}\\
&=T_{y_{0}}\big(T_{(ts')^{n-1}}T_{t}\sum_{w_{i}\in W'} a_{w_{i}}C_{w_{i}}+ v^{-L(t)}T_{(s't)^{n-1}}\sum_{w_{i}\in W'} a_{w_{i}}C_{w_{i}} \big)\\
&=\sum a_{w_{i}}T_{y_{0}.(ts')^{n-1}}C_{tw_{i}}+v^{-L(t)}T_{y_{0}}\sum a_{w_{i}}\big(T_{(s't)^{n-1}}C_{w_{i}}-  T_{(ts')^{n-1}}C_{w_{i}}\big).
\end{align*}
By induction we see that
$$\sum a_{w_{i}}T_{y_{0}}T_{(ts')^{n-1}}C_{tw_{i}}\in\cM.$$
Lemma  \ref{Tmult} implies that
$$T_{(s't)^{n-1}}C_{w}-  T_{(ts')^{n-1}}C_{w}$$ 
is an $\cA$-linear combination of terms of the form $T_{(s't)^{m}s'}C_{tw'}$ and $T_{(ts')^{m}}C_{tw'}$, for some $tw'\in U$ and $m\leq n-2$ (it is $0$ if $n=1$).  Thus it follows by induction that 
$$T_{y_{0}}\sum a_{w_{i}}\big(T_{(s't)^{n-1}}C_{w_{i}}-  T_{(ts')^{n-1}}C_{w_{i}}\big)\in\cM$$
as required.
\end{proof}
\begin{Prop}
\label{pol}
For all $u\in U$, $u_{1}<u$ and $y\in X_{u}$ we have
$$P_{u_{1},u}T_{y}T_{u_{1}}\quad \text{is an $\cA_{<0}$-linear combination of $T_{z}$.}$$
\end{Prop}
\begin{proof}
Let $u=tw'\in U$, $u_{1}<u$ and $y\in X_{u}$. One can see that we have either $u_{1}\in W'$ (then $u_{1}\leq w'$) or there exists $w\in W'$ such that $u_{1}=tw$ and $w< w'$. \\
Assume that $u_{1}\in W'$. Then $tu_{1}>u_{1}$ and we have (using (\cite[Theorem 6.6]{bible})
$$P_{u_{1},u}=P_{u_{1},tw'}=v^{-L(t)}P_{tu_{1},tw'}\in v^{-L(t)}\cA_{\leq 0}.$$
Furthermore, the degree of the polynomials occurring in the decomposition of $T_{y}T_{u_{1}}$ in the standard basis is at most $N$. Indeed, let $y'\in X_{S'}$ and $v\in W'$ be such that $y=y'v$. Then we have
\begin{align*}
T_{y}T_{u_{1}}&=T_{y'}T_{v}T_{u_{1}}\\
&=T_{y'}\sum_{u'\in W'} f_{v,u_{1},u'}T_{u'}\\
&=\sum_{u'\in W'} f_{v,u_{1},u'}T_{y'u'}
\end{align*}
and since $W'$ is bounded by $N$, the degree of $ f_{v,u_{1},u'}$ is less than or equal to $N$. Thus, since $L(t)> N$, we get the result in that case.\\
Assume that $u_{1}=tw$ ($w\in W'$). Then, since $y\in (\cR^{\{t\}}\cap X_{S'})t$, we see that $\ell(yu_{1})=\ell(y)+\ell(u_{1})$ and $T_{y}T_{u_{1}}=T_{yu_{1}}$. The result follows.
\end{proof}

We are now ready to prove Theorem \ref{tfix}. Conditions  {\bf I4} and  {\bf I5} follow respectively from Proposition \ref{mod} and \ref{pol}. Applying Theorem \ref{main} yields that
$$\{xu|\ u\in U, x\in X_{u}\}=\{w\in W|\ w=yw',\ y\in \cR^{\{t\}}\cap X_{S'}, w'\in W'\}$$
is a left ideal of $W$. 

\begin{Exa}
\label{Fs}
Let $W$ be of type $\tilde{G_{2}}$ with presentation as follows
$$W:=\sg s_{1},s_{2},s_{3}\ |\ (s_{1}s_{2})^{6}=1,(s_{2}s_{3})^{3}=1,(s_{1}s_{3})^{2}=1\sd$$
\smallskip
and let $L$ be a weight function on $W$. 
The longest element of the subgroup $W'$ generated by $s_{2},s_{3}$ is $w_{0}=s_{2}s_{3}s_{2}$ and $L(w_{0})=3L(s_{2})$. One can easily check that $3L(s_{2})$ is a bound for $W'$, thus if $L(s_{1})>3L(s_{2})$ we can apply Theorem \ref{main4}. We obtain that the following sets (which are the cells of $W'$):
\begin{align*}
&\{e\}\cup\{s_{2},s_{3}s_{2}\} \cup\{s_{3},s_{2}s_{3}\}\cup\{w_{0}\}\qquad \text{(left cells)}\\
&\{e\}\cup \{s_{2},s_{3},s_{3}s_{2},s_{2}s_{3}\}\cup \{w_{0}\}\qquad \text{(two-sided cells)}.
\end{align*}
 are left cells (resp. two-sided cells) of $W$. \\
\end{Exa}

\section{Miscellaneous}
In this section $(W,S)$ denotes an arbitrary Coxeter system and $L$ a positive weight function on $W$. We give a number of lemmas which will be needed later on. 
\begin{Lem}
\label{mu}
Let $S'\subseteq S$ be such that
\begin{enumerate}
\item for all $s'_{1},s'_{2}\in S'$, we have $L(s'_{1})=L(s'_{2})$,
\item for all $t\in S-S'$ and $s'\in S'$ we have $L(t)>L(s')$.
\end{enumerate}
Let $y,w\in W$ and $s'\in S'$ be such that $s'y<y<w<s'w$. Then if $M^{s'}_{y,w}\neq 0$, we have either $\cL(w)\subseteq \cL(y)$ or there exists $s\in S'$ such that $w=sy$, in which case $M_{y,w}^{s'}=1$.
\end{Lem}
\begin{proof}
We proceed by induction on $\ell(w)-\ell(y)$. Assume first that $\ell(w)-\ell(y)=1$. Since $s'y<y$ and $s'w>w$ one can see that there exist $s\in S$ such that $s\neq s'$ and $w=sy$. In that case we have 
$$M_{y,w}^{s'}=
\begin{cases}
0, & \mbox{if }L(s)>L(s'),\\ 
1, & \mbox{if }L(s)=L(s'). 
\end{cases}
$$
Thus if $M_{z,w}^{s'}\neq 0$ we must have $s\in S'$.\\
Assume that $\ell(w)-\ell(y)>1$ and that $\cL(w)\nsubseteq \cL(y)$. Let $s\in S$ be such that  $s\in\cL(w)$ and $s\notin\cL(y)$. We have
$$M_{y,w}^{s'}+\sum_{z;y<z<w,s'z<z}P_{y,z}M_{z,w}^{s'}-v_{s'}P_{y,w}\in\cA_{<0}.$$
Thus in order to show that $M_{y,w}^{s'}=0$ it is enough to show that
$$\sum_{z;y<z<w,s'z<z}P_{y,z}M_{z,w}^{s'}-v_{s'}P_{y,w}\in\cA_{<0}.$$
Let $z\in W$ be such that $M^{s'}_{z,w}\neq 0$. By induction we have either $M_{z,w}^{s'}=1$ or $\cL(w)\subseteq \cL(z)$. In the first case we have $P_{y,z}M_{z,w}^{s'}\in\cA_{<0}$. Assume that we are in the second case (then $s\in\cL(z)$). By (\cite[proof of Theorem 6.6]{bible}) we know that
$$P_{y,z}=v_{s}^{-1}P_{sy,z}\in \cA_{\leq 0}.$$
Furthermore the degree in $v$ of $M_{z,w}^{s'}$ is at most $L(s')-1$ (\cite[Proposition 6.4]{bible}). Since $s'\in S'$ we have $L(s)\geq L(s')$ and
$$P_{y,z}M_{z,w}^{s'}\in \cA_{<0}.$$
Similarly $v_{s'}P_{y,w}\in \cA_{<0}$ (since $\ell(w)-\ell(y)>1$). Thus if $\cL(w)\nsubseteq \cL(y)$ we must have $M_{y,w}^{s'}=0$, as required.

\end{proof}

\begin{Lem}
\label{mulc}
Let $\fB\subseteq W$ be a left ideal of $W$. Let $s\in S$ and $\fB'_{s}$ (resp. $\fB_{s}$) be the subset of $\fB$ which consists of all $w\in \fB$ such that $ws>w$ (resp. $ws<w$). Assume that there exists a left ideal $\fA$ of $W$ such that, for all $w'\in \fB'_{s}$ we have
$$C_{w'}C_{s}=C_{w's}+\sum_{z\in \fA} \cA C_{z}.$$
Then  $\fA\cup \fB_{s}\cup \fB'_{s}s$ is a left ideal of $W$. 
\end{Lem}
\begin{proof}
Let $w\in \fA\cup \fB_{s}\cup \fB'_{s}s$. Let $y\in W$ be such that $y\leq_{L} w$. We need to show that $y\in  \fA\cup \fB_{s}\cup \fB'_{s}s$.\\
If $w\in \fA$ then $y\in \fA$, since $\fA$ is a left ideal.\\
If $w\in \fB_{s}$ then $y\in\fB$. Note that since
$$y\leq_{L} w\Longrightarrow \cR(w)\subseteq \cR(y),$$
we have $s\in \cR(y)$ and $y\in \fB_{s}$. This shows that $\fB_{s}$ is a left ideal.\\
Finally, assume that $w\in \fB'_{s}s$ and let $w'=ws\in \fB'_{s}$. We may assume that there exists $t\in S$ such that $C_{y}$ appears with a non-zero coefficient in the expression of $C_{t}C_{w}$ in the Kazhdan-Lusztig basis. We have
\begin{align*}
C_{t}C_{w}&=C_{t}C_{w's}\\
&=C_{t}\big(C_{w'}C_{s}+\sum_{z\in \fA} \cA C_{z}  \big)\\
&=\big(\sum_{z\in \fB} \cA C_{z}\big) C_{s}+\sum_{z\in \fA}\cA C_{z}\\
&=\sum_{z\in \fB'_{s}s} \cA C_{z}+\sum_{z\in \fB_{s}} \cA C_{z}+\sum_{z\in \fA} \cA C_{z}
\end{align*}
Thus we see that $y\in \fA\cup \fB_{s}\cup \fB'_{s}s$ as desired. 
\end{proof}

\begin{Lem}
\label{itsc}
Let $T$ be a union of left cells which is stable by taking the inverse. Let $T=\cup \ T_{i}$ ($1\leq i\leq N$) be the decomposition of $T$ into left cells. Assume that for all $i,j\in\{1,...,N\}$ we have
\begin{equation}
T_{i}^{-1}\cap T_{j}\neq \emptyset \tag{$\ast$}
\end{equation}
Then $T$ is included in a two-sided cell.
\end{Lem}
\begin{proof}
Let $y,w\in T$ and $i,j\in \{1,...,N\}$ be such that $y\in T_{i}$ and $w\in T_{j}$. Using ($\ast$), there exist $y_{1},y_{2}\in T_{i}$ such that $y_{1}^{-1}\in T_{i}$ and $y_{2}^{-1}\in T_{j}$. We have
$$y\sim_{L} y_{1}\sim_{L} y_{2}\quad\Longrightarrow\quad y\sim_{L} y_{1}^{-1}\sim_{R} y_{2}^{-1}\sim_{L} w$$
as required.
\end{proof}

\section{Decomposition of $\tilde{G}_{2}$ in the asymptotic case}

Let $W$ be an affine Weyl group of type $\tilde{G_{2}}$ with diagram and weight function given by
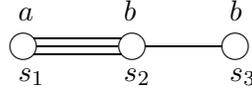
\begin{center}
\begin{picture}(150,32)
\put( 40, 10){\circle{10}}
\put( 44,  7){\line(1,0){33}}
\put( 45,  10){\line(1,0){30.5}}
\put( 44, 13){\line(1,0){33}}
\put( 81, 10){\circle{10}}
\put( 86, 10){\line(1,0){29}}
\put(120, 10){\circle{10}}
\put( 38, 20){$a$}
\put( 78, 20){$b$}
\put(118, 20){$b$}
\put( 38, -3){$s_{1}$}
\put( 78, -3){$s_{2}$}
\put(118,-3){$s_{3}$}
\end{picture}
\end{center}
where $a,b$ are positive integers.

The aim of this section is to find the decomposition of $W$ into left cells and two-sided cells for any weight function $L$ such that $a/b>4$. Furthermore we will determine the partial left (resp. two-sided) order on the left (resp. two-sided) cells (see Section \ref{leftorder}). We fix such a weight function $L$. Throughout this section, we keep this setting.

In Figure 2, we present a partition of $W$ using the geometric realization as described in Section \ref{geom}, where the pieces are formed by the alcoves lying in the same connected component after removing the thick lines. We have
\begin{Th}
The partition of $W$ described in Figure 2 coincides with the partition of $W$ into left cells.
\end{Th}
Using the same methods as in \cite[Section 6]{jeju}, one can show that each of the pieces is included in a left cell (with respect to $L$). Thus in order to prove that each of the pieces is a left cell it is enough to show that each of them is included in a union of left cells.\\

We now consider the union of all subsets of $W$ whose name contains a fixed capital letter; we denote this union by that capital letter. For instance
$$A=(\overset{6}{\underset{i=1}{\cup}} A_{i})\ \bigcup\  (\overset{6}{\underset{i=1}{\cup}} A'_{i}).$$ 
We have
\begin{Th}
The decomposition of $W$ into two-sided cells is as follows
$$W=A\cup B\cup C\cup D\cup E\cup F\cup \{e\}.$$
\end{Th}

The proof of these theorems will be given in the next sections. For a start, we already know that (see \cite[\S 4]{jeju2})
\begin{itemize}
\item $A$ is a two sided cell;
\item $A_{i}$ and $A'_{i}$ are left cells for all $1\leq i \leq 6$;
\item $A_{i}$ and $A'_{i}$ are left ideals for all $1\leq i \leq 6$.
\end{itemize}

\begin{Rem}
\label{comput}
In this section we need to compute some Kazhdan-Lusztig polynomials $P_{x,y}$ ($x,y\in W$) for a whole class of weight functions. Methods for dealing with this problem are presented in \cite[Proposition 3.2 and \S6]{jeju}. In particular, this involved some computations with GAP (\cite{GAP}).
\end{Rem}


%

%


We now recall some notation. For any subset $J\subseteq \{s_{1},s_{2},s_{3}\}$, let
\begin{enumerate}
\item $\cR^{J}:=\{w\in W\ |\ J \subseteq \cR(w)\}$;
\item $W_{J}$ be the subgroup of $W$ generated by $J$;
\item $X_{J}:=\{w\in W|\ w \text{ has minimal length in $wW_{J}$}\}$.
\end{enumerate}
$\ $\\
We refer to \cite{webpage} for details in the computations. 

$\ $\\
\begin{center}
\begin{pspicture}(-6.2,6.2)(6.2,-6.2)
\psset{unit=0.90cm}
\SpecialCoor
\psline(0,-6)(0,6)
\multido{\n=1+1}{7}{
\psline(\n;30)(\n;330)}
\rput(1;30){\psline(0,0)(0,5.5)}
\rput(1;330){\psline(0,0)(0,-5.5)}
\rput(2;30){\psline(0,0)(0,5)}
\rput(2;330){\psline(0,0)(0,-5)}
\rput(3;30){\psline(0,0)(0,4.5)}
\rput(3;330){\psline(0,0)(0,-4.5)}
\rput(4;30){\psline(0,0)(0,4)}
\rput(4;330){\psline(0,0)(0,-4)}
\rput(5;30){\psline(0,0)(0,3.5)}
\rput(5;330){\psline(0,0)(0,-3.5)}
\rput(6;30){\psline(0,0)(0,3)}
\rput(6;330){\psline(0,0)(0,-3)}
\rput(7;30){\psline(0,0)(0,2.5)}
\rput(7;330){\psline(0,0)(0,-2.5)}
\multido{\n=1+1}{7}{
\psline(\n;150)(\n;210)}
\rput(1;150){\psline(0,0)(0,5.5)}
\rput(1;210){\psline(0,0)(0,-5.5)}
\rput(2;150){\psline(0,0)(0,5)}
\rput(2;210){\psline(0,0)(0,-5)}
\rput(3;150){\psline(0,0)(0,4.5)}
\rput(3;210){\psline(0,0)(0,-4.5)}
\rput(4;150){\psline(0,0)(0,4)}
\rput(4;210){\psline(0,0)(0,-4)}
\rput(5;150){\psline(0,0)(0,3.5)}
\rput(5;210){\psline(0,0)(0,-3.5)}
\rput(6;150){\psline(0,0)(0,3)}
\rput(6;210){\psline(0,0)(0,-3)}
\rput(7;150){\psline(0,0)(0,2.5)}
\rput(7;210){\psline(0,0)(0,-2.5)}
\multido{\n=1.5+1.5}{4}{
\psline(-6.062,\n)(6.062,\n)}
\multido{\n=0+1.5}{5}{
\psline(-6.062,-\n)(6.062,-\n)}
\psline(0;0)(7;30)

\rput(0,1){\psline(0;0)(7;30)}
\rput(0,1){\psline(0;0)(7;210)}
\rput(0,2){\psline(0;0)(7;30)}
\rput(0,2){\psline(0;0)(7;210)}
\rput(0,3){\psline(0;0)(6;30)}
\rput(0,3){\psline(0;0)(7;210)}
\rput(0,4){\psline(0;0)(4;30)}
\rput(0,4){\psline(0;0)(7;210)}
\rput(0,5){\psline(0;0)(2;30)}
\rput(0,5){\psline(0;0)(7;210)}
\rput(0,6){\psline(0;0)(7;210)}
\rput(-1.732,6){\psline(0;0)(5;210)}
\rput(-3.464,6){\psline(0;0)(3;210)}
\rput(-5.196,6){\psline(0;0)(1;210)}

\rput(0,-1){\psline(0;0)(7;30)}
\rput(0,-1){\psline(0;0)(7;210)}
\rput(0,-2){\psline(0;0)(7;30)}
\rput(0,-2){\psline(0;0)(7;210)}
\rput(0,-3){\psline(0;0)(7;30)}
\rput(0,-3){\psline(0;0)(6;210)}
\rput(0,-4){\psline(0;0)(7;30)}
\rput(0,-4){\psline(0;0)(4;210)}
\rput(0,-5){\psline(0;0)(7;30)}
\rput(0,-5){\psline(0;0)(2;210)}
\rput(0,-6){\psline(0;0)(7;30)}
\rput(1.732,-6){\psline(0;0)(5;30)}
\rput(3.464,-6){\psline(0;0)(3;30)}
\rput(5.196,-6){\psline(0;0)(1;30)}

\psline(0;0)(6.928;60)
\rput(1.732,0){\psline(0;0)(6.928;60)}
\rput(3.464,0){\psline(0;0)(5.196;60)}
\rput(5.196,0){\psline(0;0)(1.732;60)}
\rput(1.732,0){\psline(0;0)(6.928;240)}
\rput(3.464,0){\psline(0;0)(6.928;240)}
\rput(5.196,0){\psline(0;0)(6.928;240)}
\rput(6.062,-1.5){\psline(0;0)(5.196;240)}
\rput(6.062,-4.5){\psline(0;0)(1.732;240)}
\rput(-1.732,0){\psline(0;0)(6.928;60)}
\rput(-3.464,0){\psline(0;0)(6.928;60)}
\rput(-5.196,0){\psline(0;0)(6.928;60)}
\rput(-1.732,0){\psline(0;0)(6.928;240)}
\rput(-3.464,0){\psline(0;0)(5.196;240)}
\rput(-5.196,0){\psline(0;0)(1.732;240)}
\rput(-6.062,1.5){\psline(0;0)(5.196;60)}
\rput(-6.062,4.5){\psline(0;0)(1.732;60)}

\psline(0;0)(6.928;120)
\psline(0;0)(6.928;300)
\rput(-1.732,0){\psline(0;0)(6.928;120)}
\rput(-3.464,0){\psline(0;0)(5.196;120)}
\rput(-5.196,0){\psline(0;0)(1.732;120)}
\rput(-1.732,0){\psline(0;0)(6.928;300)}
\rput(-3.464,0){\psline(0;0)(6.928;300)}
\rput(-5.196,0){\psline(0;0)(6.928;300)}
\rput(-6.062,-1.5){\psline(0;0)(5.196;300)}
\rput(-6.062,-4.5){\psline(0;0)(1.732;300)}
\rput(1.732,0){\psline(0;0)(6.928;300)}
\rput(3.464,0){\psline(0;0)(5.196;300)}
\rput(5.196,0){\psline(0;0)(1.732;300)}
\rput(1.732,0){\psline(0;0)(6.928;120)}
\rput(3.464,0){\psline(0;0)(6.928;120)}
\rput(5.196,0){\psline(0;0)(6.928;120)}
\rput(6.062,1.5){\psline(0;0)(5.196;120)}
\rput(6.062,4.5){\psline(0;0)(1.732;120)}

\rput(0,1){\psline(0;0)(7;150)}
\rput(0,1){\psline(0;0)(7;330)}
\rput(0,2){\psline(0;0)(7;150)}
\rput(0,2){\psline(0;0)(7;330)}
\rput(0,3){\psline(0;0)(6;150)}
\rput(0,3){\psline(0;0)(7;330)}
\rput(0,4){\psline(0;0)(4;150)}
\rput(0,4){\psline(0;0)(7;330)}
\rput(0,5){\psline(0;0)(7;330)}
\rput(0,6){\psline(0;0)(7;330)}
\rput(-1.732,6){\psline(0;0)(2;330)}

\rput(0,-1){\psline(0;0)(7;150)}
\rput(0,-1){\psline(0;0)(7;330)}
\rput(0,-2){\psline(0;0)(7;150)}
\rput(0,-2){\psline(0;0)(7;330)}
\rput(0,-3){\psline(0;0)(7;150)}
\rput(0,-3){\psline(0;0)(6;330)}
\rput(0,-4){\psline(0;0)(7;150)}
\rput(0,-4){\psline(0;0)(4;330)}
\rput(0,-5){\psline(0;0)(7;150)}
\rput(0,-5){\psline(0;0)(2;330)}
\rput(0,-6){\psline(0;0)(7;150)}

\rput(6.062,3.5){\psline(0;0)(5;150)}
\rput(6.062,4.5){\psline(0;0)(3;150)}
\rput(6.062,5.5){\psline(0;0)(1;150)}

\rput(-6.062,-3.5){\psline(0;0)(5;330)}
\rput(-6.062,-4.5){\psline(0;0)(3;330)}
\rput(-6.062,-5.5){\psline(0;0)(1;330)}

\psline(0;0)(7;330)
\psline(0;0)(7;150)
\psline(0;0)(7;210)
\psline(0;0)(6.928;240)

\pspolygon[linewidth=0.7mm](0,0)(0,1)(0.433,0.75)
\pspolygon[linewidth=0.7mm](0,0)(0,1)(-0.866,1.5)
\pspolygon[linewidth=0.7mm](0,0)(-0.433,0.75)(-1.732,0)(-0.866,0)(-0.866,-1.5)
\pspolygon[linewidth=0.7mm](0,0)(0.433,0.75)(1.732,0)(0.866,0)(0.866,-1.5)
\pspolygon[linewidth=0.7mm](0,1)(0,1.5)(-0.866,1.5)
\pspolygon[linewidth=0.7mm](0,1)(0,3)(0.433,2.25)(1.732,3)(0.433,0.75)
\psline[linewidth=0.7mm](0,0)(0,-6)
\psline[linewidth=0.7mm](0;0)(6.928;240)
\rput(0.866,-1.5){\psline[linewidth=0.7mm](0;0)(4.5;270)}
\rput(0.866,-1.5){\psline[linewidth=0.7mm](0;0)(5.196;300)}
\rput(1.732,0){\psline[linewidth=0.7mm](0;0)(6.928;300)}
\rput(2.598,-1.5){\psline[linewidth=0.7mm](0;0)(4;330)}
\rput(1.732,0){\psline[linewidth=0.7mm](0;0)(5;330)}
\rput(1.732,0){\psline[linewidth=0.7mm](0;0)(4.33;0)}
\rput(0.866,1.5){\psline[linewidth=0.7mm](0;0)(5.196;0)}
\rput(0.866,1.5){\psline[linewidth=0.7mm](0;0)(6;30)}
\rput(2.598,1.5){\psline[linewidth=0.7mm](0;0)(4;30)}
\rput(1.732,3){\psline[linewidth=0.7mm](0;0)(3.464;60)}
\rput(0,3){\psline[linewidth=0.7mm](0;0)(3.464;60)}
\rput(0.866,4.5){\psline[linewidth=0.7mm](0;0)(1.5;90)}
\rput(0,3){\psline[linewidth=0.7mm](0;0)(3;90)}
\rput(0,3){\psline[linewidth=0.7mm](0;0)(3.464;120)}
\rput(-0.866,1.5){\psline[linewidth=0.7mm](0;0)(5.196;120)}
\rput(-0.866,1.5){\psline[linewidth=0.7mm](0;0)(6;150)}
\rput(-0.866,1.5){\psline[linewidth=0.7mm](0;0)(5.196;180)}
\rput(-2.598,1.5){\psline[linewidth=0.7mm](0;0)(4;150)}
\rput(-1.732,0){\psline[linewidth=0.7mm](0;0)(4.33;180)}
\rput(-1.732,0){\psline[linewidth=0.7mm](0;0)(5;210)}
\rput(-1.732,0){\psline[linewidth=0.7mm](0;0)(6.928;240)}
\rput(-1.732,0){\psline[linewidth=0.7mm](0;0)(6.928;240)}
\rput(-2.598,-1.5){\psline[linewidth=0.7mm](0;0)(4;210)}
\rput(0,1.5){\psline[linewidth=0.7mm](0;1.5)(0.866;1.5)}

\rput(0.17,0.68){{\tiny $A_{0}$}}
\rput(-1.5,-4.33){{\small $A_{1}$}}
\rput(-4.97,-4.3){{\small $A'_{6}$}}
\rput(-4.97,-1.3){{\small $A_{6}$}}
\rput(-4.97,1.7){{\small $A'_{5}$}}
\rput(-2.372,3.2){{\small $A_{5}$}}
\rput(-0.23,4.67){{\small $A_{4}$}}
\rput(1.1,5.85){{\small $A'_{4}$}}
\rput(2.824,3.2){{\small $A_{3}$}}
\rput(3.69,1.7){{\small $A'_{3}$}}
\rput(5.452,-1.31){{\small $A_{2}$}}
\rput(5.452,-4.31){{\small $A'_{2}$}}
\rput(1.988,-4.31){{\small $A'_{1}$}}
\psline[linewidth=0.5mm]{->}(0.433,5.75)(0.433,6.75)
\rput(0.433,7){$C_{4}$}
\rput(5.812,3.8){\psline[linewidth=0.5mm]{->}(0;0)(1;30)}
\rput(7,4.4){$C_{3}$}
\rput(5.812,-2.8){\psline[linewidth=0.5mm]{->}(0;0)(1;330)}
\rput(7,-3.4){$C_{2}$}
\psline[linewidth=0.5mm]{->}(0.433,-5.75)(0.433,-6.75)
\rput(0.433,-7){$C_{1}$}
\rput(-5.812,-2.8){\psline[linewidth=0.5mm]{->}(0;0)(1;210)}
\rput(-7,-3.4){$C_{6}$}
\rput(-5.812,3.8){\psline[linewidth=0.5mm]{->}(0;0)(1;150)}
\rput(-7,4.4){$C_{5}$}
\rput(2.4,5.75){\psline[linewidth=0.5mm]{->}(0;0)(1;60)}
\rput(3.05,6.83){$B_{3}$}
\rput(5.812,0.75){\psline[linewidth=0.5mm]{->}(0;0)(1;0)}
\rput(7.1,0.75){$B_{2}$}
\rput(4.132,-5.75){\psline[linewidth=0.5mm]{->}(0;0)(1;300)}
\rput(4.837,-6.83){$B_{1}$}
\rput(-4.132,-5.75){\psline[linewidth=0.5mm]{->}(0;0)(1;240)}
\rput(-4.837,-6.83){$B_{6}$}
\rput(-5.812,0.75){\psline[linewidth=0.5mm]{->}(0;0)(1;180)}
\rput(-7.1,0.75){$B_{5}$}
\rput(-2.4,5.75){\psline[linewidth=0.5mm]{->}(0;0)(1;120)}
\rput(-3.05,6.83){$B_{4}$}

\rput(0.22,1.33){\tiny $E_{1}$}

\rput(-0.32,1.01){\tiny $E_{2}$}

\rput(0.33,1.98){\tiny $D_{3}$}

\rput(-0.65,0.19){\tiny $D_{2}$}

\rput(0.65,0.17){\tiny $D_{1}$}

\rput(-.25,1.35){\tiny $F$}

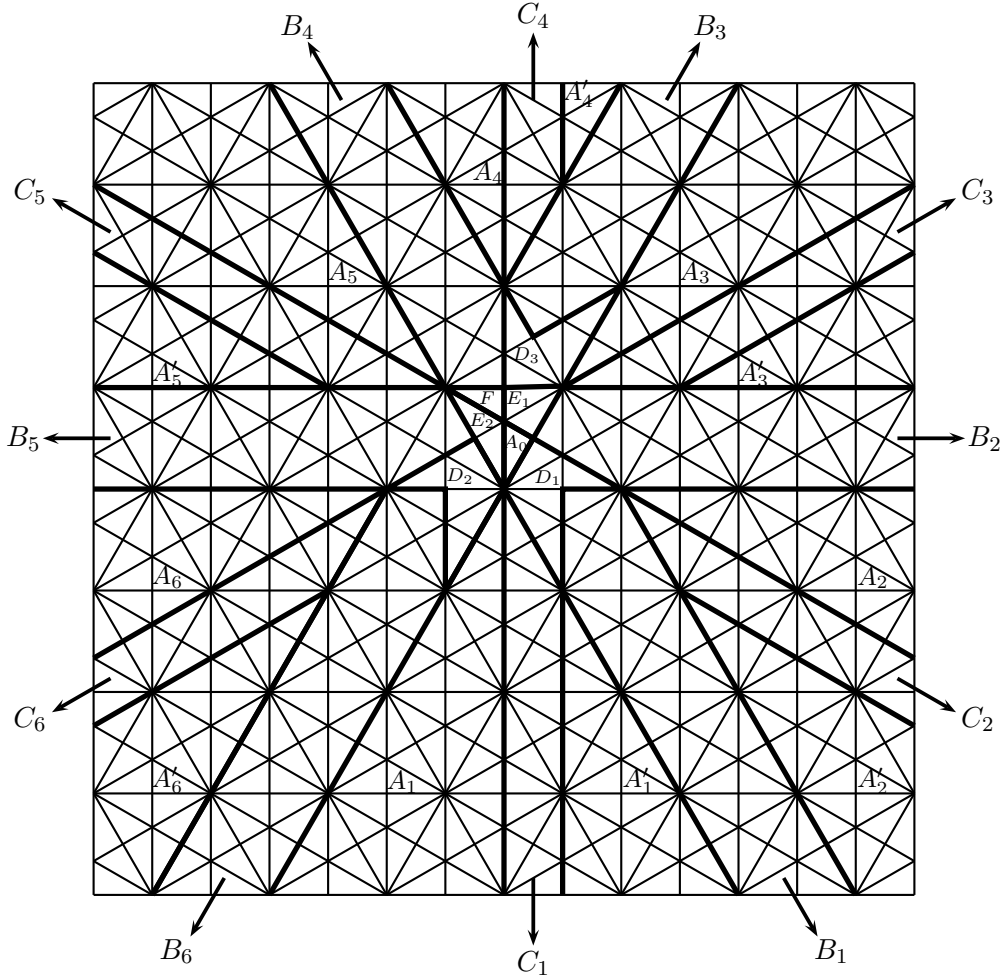
\rput(0,-7.8){Figure 2. Decomposition of $\tilde{G}_{2}$ into left cells in the case $a>4b$}

\end{pspicture}
\end{center}

$\ $\\

\subsection{The sets $C_{i}$}
\label{setC}
In this section we want to prove that $C_{i}$ (for all $1\leq i\leq 6$) is a left cell and that $C=\cup C_{i}$ is a two-sided cell. 

For $1\leq i\leq 6$, let
\begin{enumerate}
\item $u_{i}\in C_{i}$ be the element of minimal length in $C_{i}$;
\item $v_{i}\in A_{i}$ be the element of minimal length in $A_{i}$;
\item $v'_{i}\in A'_{i}$ be the element of minimal length in $A'_{i}$.
\end{enumerate}
For instance, we have 
\begin{align*}
u_{1}&=s_{1}s_{2}s_{1}s_{2}s_{1};\\
v_{1}&=s_{1}s_{2}s_{1}s_{2}s_{1}s_{2};\\
v'_{1}&=s_{2}s_{1}s_{2}s_{1}s_{2}s_{3}s_{1}s_{2}s_{1}s_{2}s_{1}.
\end{align*}
We set $U:=\{u_{i},v_{i},v'_{i}\ |\ 1\leq i\leq 6\}$, $X_{v_{i}}=X_{v'_{i}}=X_{\{s_{1},s_{2}\}}$ and
$$X_{u_{i}}=\{z\in W\ |\ zu_{i}\in C_{i}\}$$
for all $1\leq i\leq 6$. We want to apply Corollary \ref{mainC}. One can check that conditions ${\bf I1}$--${\bf I3}$ of Theorem \ref{main} hold. We now have a look at condition {\bf I4}.
\begin{Lem} 
\label{Cs}
The submodule
$$\cM:=\sg T_{x}C_{u}\ |\ u\in U, x\in X_{u}\sd_{\cA}\subseteq \cH$$
is a left ideal.
\end{Lem}
\begin{proof}
In \cite[Lemma 5.2]{jeju2}, it has been shown that 
$$\sg T_{x}C_{v_{i}}\ |\  x\in X_{\{s_{1},s_{2}\}}\sd_{\cA}\quad\text{and}\quad \sg T_{x}C_{v'_{i}}\ |\  x\in X_{\{s_{1},s_{2}\}}\sd_{\cA}$$
are left ideals of $\cH$, for all $1\leq i\leq 6$. Thus, in order to show that $\cM$ is a left ideal of $\cH$, it is enough to prove that $T_{x}C_{u_{i}}\in\cM$ for all $1\leq i\leq 6$ and all $x\in W$. We proceed by induction on $\ell(x)$. If $\ell(x)=0$ it's clear. Assume that $\ell(x)>0$. We may assume that $x\notin X_{u_{i}}$. Then, one can see that we have either $x=x_{0}s_{2}$ (and $\ell(x)=\ell(x_{0})+1$) or $x=x_{1}s_{2}s_{1}s_{2}s_{1}s_{2}s_{3}$ (and $\ell(x)=\ell(x_{1})+6$). 
Now, doing some explicit computations, one can show that $T_{s_{2}}C_{u_{i}}$ is an $\cA$-linear combination of $C_{u}$ with $u\in U$. For example, we have 
$$T_{s_{2}}C_{u_{1}}=C_{v_{1}}-v^{-L(s_{2})}C_{u_{1}}$$
and
$$T_{s_{2}}C_{u_{5}}=C_{v_{5}}-v^{-L(s_{2})}C_{u_{5}}+C_{v_{1}}.$$
Similarly, one can show that $T_{s_{2}s_{1}s_{2}s_{1}s_{2}s_{3}}C_{u_{i}}$ is an $\cA$-linear combination of 
terms of the form $T_{z}C_{u}$ where $u\in U$, $z\in X_{u}$ and $\ell(z)<\ell(s_{2}s_{1}s_{2}s_{1}s_{2}s_{3})$. For example we have
\begin{align*}
T_{s_{2}s_{1}s_{2}s_{1}s_{2}s_{3}}C_{u_{1}}&=C_{v'_{1}}+\cA T_{s_{1}s_{2}s_{1}s_{2}s_{3}}C_{u_{1}}+\cA T_{s_{2}s_{1}s_{2}s_{3}}C_{u_{1}}+\cA T_{s_{1}s_{2}s_{3}}C_{u_{1}}\\
&\ \ \ \ \ +\cA T_{s_{2}s_{3}}C_{u_{1}}+\cA T_{s_{3}}C_{u_{1}}+\cA C_{u_{1}}+\cA C_{v_{1}}.\\
\end{align*}
Thus by induction, we obtain that $T_{x}C_{u_{i}}\in \cM$ as required. 
\end{proof}

We now have a look at condition {\bf I5}. Let $u\in U$, $u'<u$ and $y\in X_{u}$. We need to show that
$$P_{u',u}T_{y}T_{u'} \text{ is an $\cA_{<0}$-linear combination of $T_{z}$}.$$
For $u=v_{i}$ or $u=v'_{i}$, it has been proved in \cite[Lemma 5.1]{jeju2}. In order to prove it for $u=u_{i}$ we proceed as follows. We determine an upper bound for the degree of the polynomials occurring in the expression of $T_{y}T_{u'}$ (where $y\in C_{i}$, $u'<u_{i}$) in the standard basis using either \cite[Theorem 2.1]{jeju2} or explicit computations. Then we compute the polynomials $P_{u',u}$ (see Remark \ref{comput}) and we can check that the condition is satisfied for all weight functions such that $L(s_{1})>4L(s_{2})$. 

We can now apply Corollary \ref{mainC}. We need to find the equivalence classes on $U$ with respect to $\preceq$.
Using the fact that $\sg T_{x}C_{v_{i}}\ |\  x\in X_{\{s_{1},s_{2}\}}\sd_{\cA}$ and $\sg T_{x}C_{v'_{i}}\ |\  x\in X_{\{s_{1},s_{2}\}}\sd_{\cA}$ are left ideals of $\cH$ for all $1\leq i\leq 6$ and  the relations computed in the previous proof, one can check that
$$\{\{u_{i}\}\{v_{i}\}, \{v'_{i}\}\ |\ 1\leq i\leq 6\}$$
is the decomposition of $U$ into equivalence classes. Hence by Corollary \ref{mainC}, the set $X_{u_{i}}u_{i}=C_{i}$ is a union of left cells for all $1\leq i\leq 6$. Since $C_{i}$ is included in a left cell, we obtain that each of the $C_{i}$'s is a left cell. \\
More precisely, if $L$ is a weight function such that $a/b>4$,  the following sets are left ideals of $W$
$$C_{i}\cup A_{i}\cup A'_{i}\tx{ for $i=1,2,3,6$}$$
$$C_{4}\cup A_{4}\cup A'_{4}\cup A_{2},$$
$$C_{5}\cup A_{5}\cup A'_{5}\cup A_{1}.$$

\begin{Prop}
\label{2sC}
The set $C$ is a two-sided cell.  
\end{Prop}
\begin{proof}
Applying Theorem \ref{main} to the set $U$ yields that $A\cup C$ is a left ideal of $W$. One can check that $A\cup C$ is stable by taking the inverse, thus it is a two-sided ideal and $A\cup C$ is a union of two-sided cells. Since $A$ is a two sided cell (see \cite{jeju2} and the references there), we see that $C$ is a union of two-sided cells. Now one can check that $C=\cup C_{i}$ satisfy the requirement of Lemma \ref{itsc} thus $C$ is included in a two-sided cell. It follows that $C$ is a two-sided cell.
\end{proof}

\subsection{The sets $B_{i}$} We want to prove that $B_{i}$ (for all $1\leq i\leq 6$) is a left cell. To this end, since $B_{i}$ is included in a left cell, it is enough to show that $B_{i}$ is a union of left cells. We also show that $B$ is a two-sided cell.

\begin{Cl}
\label{B1}
The set $B_{1}$ is a left cell.
\end{Cl}
\begin{proof}
Set $u=s_{1}s_{3}s_{2}s_{1}$ and 
$$X_{u_{1}}=\{z\in W\ |\ zs_{1}s_{3}s_{2}s_{1}\in B_{1}\}.$$
Recall that
\begin{align*}
u_{1}&=s_{1}s_{2}s_{1}s_{2}s_{1},\\
v_{1}&=s_{1}s_{2}s_{1}s_{2}s_{1}s_{2},\\
v'_{1}&=s_{1}s_{2}s_{1}s_{2}s_{1}s_{2}s_{3}s_{2}s_{1}s_{2}s_{1},\\
u_{2}&=s_{1}s_{2}s_{1}s_{2}s_{1}s_{3}s_{2}s_{1},\\
v_{2}&=s_{2}s_{1}s_{2}s_{1}s_{2}s_{1}s_{3}s_{2}s_{1},\\
v'_{2}&=s_{2}s_{1}s_{2}s_{1}s_{2}s_{3}s_{1}s_{2}s_{1}s_{2}s_{1}s_{3}s_{2}s_{1}\\
v_{3}&=s_{2}s_{1}s_{2}s_{1}s_{2}s_{1}s_{3},
\end{align*}
and
\begin{align*}
X_{u_{i}}&=\{z\in W\ |\ zu_{i}\in C_{i}\},\\
X_{v_{i}}&=X_{v'_{i}}=X_{\{s_{1},s_{2}\}}\quad\text{for $1\leq i\leq 6$}.
\end{align*}
Using similar arguments as in Lemma \ref{Cs} and the results in Section \ref{setC}, one can check that we can apply Theorem \ref{main} to $U:=\{u,u_{1},v_{1},v'_{1},u_{2},v_{2},v'_{2},v_{3}\}$. We obtain that 
$$\{xu\ |\ u\in U, x\in X_{u}\}=A_{2}\cup A'_{2}\cup C_{2}\cup B_{1}\cup A_{1}\cup A'_{1}\cup C_{1}\cup A_{3}$$
is a left ideal. Since $A_{1}$, $A'_{i}$ and $C_{i}$ are left cells for all $1\leq i\leq 6$ it follows that $B_{1}$ is a left cell.
\end{proof}

\begin{Cl}
\label{B2}
$B_{2}$ is a left cell.
\end{Cl}
\begin{proof}
The set $\cR^{\{s_{1},s_{3}\}}$ is a left ideal of $W$ (see Example \ref{Lid}). Since we have
$$\cR^{\{s_{1},s_{3}\}}=B_{2}\cup A_{3}\cup A'_{3}\cup A_{2}\cup C_{3}$$
one can see that $B_{2}$ is a left cell. 
\end{proof}

\begin{Cl}
\label{B3}
The set $B_{3}$ is a left cell.
\end{Cl}
\begin{proof}
Let $v=s_{1}s_{3}s_{2}s_{1}s_{2}s_{3}$ and 
$$X_{v}:=\{z\in W|zv\in B_{3}\}\quad\quad\quad Y_{v}:=\{y\in X_{v}| \ell(ys_{2}s_{1}s_{2})=\ell(y)-3\}.$$
 We want to apply Theorem \ref{main} to the set $U=\{v,u_{4},v_{4},v'_{4},v_{3},v_{2},v_{5}\}$ and the corresponding $X_{u}$. Arguing as in Section \ref{setC}, one can show that conditions ${\bf I1}$--${\bf I4}$ hold. 
However, condition  {\bf I5} does not hold if (and only if) $v'=s_{1}s_{2}s_{1}s_{2}s_{3}<v$ and $y\in Y_{v}$. Indeed, let $y\in Y_{v}$ and $y_{0}=ys_{2}s_{1}s_{2}$, then we have $P_{v',v}=v^{-L(s_{3})}$ and
\begin{align*}
T_{y_{0}}T_{s_{2}s_{1}s_{2}}T_{v'}&=T_{y_{0}}\big( T_{s_{1}s_{2}s_{1}s_{2}s_{1}s_{3}}+(v^{L(s_{2})}-v^{-L(s_{2})})T_{s_{1}s_{2}s_{1}s_{2}s_{1}s_{2}s_{3}}\big)\\
&=T_{y_{0}s_{1}s_{2}s_{1}s_{2}s_{1}s_{3}}+(v^{L(s_{2})}-v^{-L(s_{2})})T_{y_{0}s_{1}s_{2}s_{1}s_{2}s_{1}s_{2}s_{3}}
\end{align*}
However, we can certainly construct the elements $\tilde{C}_{xu}$ (see the proof of Proposition \ref{tilde}) such that
$$\tilde{C}_{xu}=\overline{\tilde{C}_{xu}}\quad\text{ for all $u\in U$ and $x\in X_{u}$.}$$
Using Section \ref{setC} and doing some computations, one can check that
\begin{enumerate}
\item $\tilde{C}_{xu}=C_{xu}$ for all $u\in U-\{v\}$ and $x\in X_{u}$. 
\item $\tilde{C}_{yv}=C_{yv}$ if $y\in X_{v}-Y_{v}$.
\end{enumerate}
Let $y\in Y_{v}$ and $y_{0}=ys_{2}s_{1}s_{2}$. We have
$$
\begin{array}{rllll}
\tilde{C_{yv}}&=&T_{y}C_{v}+\ind{xu\sqs yv}{\ind{u\in U, x\in X_{u}}{\sum}}p^{*}_{xu,yv}T_{x}C_{u}&\\
&=&T_{y}C_{v}+\ind{x< y,x\in X_{v}}{\sum}p^{*}_{xv,yv}T_{x}C_{v}+ \ind{u\neq v}{\ind{u\in U, x\in X_{u}}{\sum}}p^{*}_{xu,yv}T_{x}C_{u}&\\
&=&T_{y}C_{v}+\ind{x< y, x\in X_{v}}{\sum}p^{*}_{xv,yv}T_{x}C_{v} \qquad\qquad\qquad& \text{mod $\cH_{<0}$}\\
&=&T_{y}C_{v} \qquad\qquad\qquad\qquad\qquad\qquad\  &\text{mod $\cH_{<0}$}\\
&=&T_{y}T_{v}+T_{y}(P_{v',v}T_{v'}) \qquad\qquad\qquad &\text{mod $\cH_{<0}$}\\
&=&T_{y}T_{v}+T_{y_{0}}T_{s_{1}s_{2}s_{1}s_{2}s_{1}s_{2}s_{3}} & \text{mod $\cH_{<0}$}\\
&=&T_{yv}+T_{y_{0}s_{1}s_{2}s_{1}s_{2}s_{1}s_{2}s_{3}}& \text{mod $\cH_{<0}$}\\
\end{array}
$$
where $\cH_{<0}=\oplus_{w\in W} \cA_{<0}T_{w}$.
Thus since $\tilde{C_{yv}}$ is stable under the involution $\ \bar{ }\ $, it follows that
$$\tilde{C_{yv}}=C_{yv}+C_{y_{0}s_{1}s_{2}s_{1}s_{2}s_{1}s_{2}s_{3}}.$$
Furthermore, since $y_{0}s_{1}s_{2}s_{1}s_{2}s_{1}s_{2}s_{3}\in A_{3}$ we obtain that
$$\sg T_{x}C_{u}| u\in U, x\in X_{u} \sd_{\cA}=\sg C_{xu} | u\in U, x\in X_{u}\sd_{\cA}$$
is a left ideal of $\cH$. We get that
$$B_{3}\cup C_{4}\cup A_{4}\cup A'_{4}\cup A_{3}\cup A_{2}\cup A_{5}$$
is a left ideal of $W$. It follows that $B_{3}$ is a left cell.
\end{proof}

\begin{Cl}
\label{B4}
The set $B_{4}$ is a left cell.
\end{Cl}
\begin{proof}
The set $\cR^{\{s_{2},s_{3}\}}$ is a left ideal of $W$. Furthermore, we have
$$\cR^{\{s_{2},s_{3}\}}=\{s_{2}s_{3}s_{2}\}\cup B_{4}\cup A_{4}\cup A_{5},$$
it follows that $B_{4}$ is a left cell.
\end{proof}
\begin{Rem}
We have seen in Example \ref{Fs} that $W-W_{\{s_{2},s_{3}\}}$ is a left ideal. Thus
$$\cR^{\{s_{2},s_{3}\}}\cap (W-W_{\{s_{2},s_{3}\}})=B_{4}\cup A_{4}\cup A_{5}$$
is a left ideal of $W$.
\end{Rem}

\begin{Cl}
\label{B5}
$B_{5}$ is a left cell.
\end{Cl}
\begin{proof}
Let $w\in \cR^{\{s_{1},s_{3}\}}$ and let $w'=ws_{1}s_{3}$. We have $ws_{2}>w$ and
$$C_{w}C_{s_{2}}=C_{ws_{2}}+\sum_{z\in W, zs_{2}<z} \mu_{z,w}^{s_{2},r}C_{z}.$$
Applying Lemma \ref{mu} (in its right version), if $M_{z,w}^{s_{2},r}\neq 0$ then we have either $\{s_{1},s_{2},s_{3}\}\subseteq \cR(z)$ which is impossible or there exists $w''\in W$ such that
$$w=w''s_{2}s_{3}\quad\text{and}\quad z=w''s_{2}.$$
Since $w=w''s_{2}s_{3}=w's_{1}s_{3}$ we must have $w\in A_{3}$, which, in turn, implies that $z\in A_{1}$ (recall that $A_{1}$ is a left ideal). Thus applying Lemma \ref{mulc} to $\fA=A_{1}$ and $\fB=\cR^{\{s_{1},s_{3}\}}$ yields that 
$$\cR^{\{s_{1},s_{3}\}}s_{2} \cup A_{1}=A_{1}\cup A_{5}\cup A'_{5}\cup A_{6}\cup C_{5}\cup B_{5}$$ 
is a left ideal of $W$. In particular $B_{5}$ is a left cell.
\end{proof}

\begin{Cl}
\label{B6}
The set $B_{6}$ is a left cell.
\end{Cl}
\begin{proof}
Applying Lemma \ref{mulc} (in a similar way as in \ref{B5}) to 
$$\fB=A_{2}\cup A'_{2}\cup C_{2}\cup B_{1}\cup A_{1}\cup A'_{1}\cup C_{1}\cup A_{3}$$
and $\fA=A_{1}$ we obtain that
$$A_{1}\cup A'_{1}\cup C_{1}\cup A_{6}\cup A'_{6}\cup C_{6}\cup A_{5}\cup B_{6}$$
is a left ideal. Thus $B_{6}$ is a left cell. In fact, since the elements of $C_{1}$ and $A'_{1}$ do not contain $s_{1}$ in their right descent set, we see that
$$A_{1}\cup A_{6}\cup A'_{6}\cup C_{6}\cup A_{5}\cup B_{6}$$
is a left ideal of $W$. 
\end{proof}

\begin{Prop}
The set $B=\cup B_{i}$ is a two-sided cell.
\end{Prop}
\begin{proof}
By the previous proofs, we see that $A\cup C\cup B$ is left ideal of $W$. Arguing as in the proof of Proposition \ref{2sC}, we obtain that $B$ is a two-sided cell.
\end{proof}


\subsection{Finite cells}

We already know that $E_{1}$, $E_{2}$, $F$ and $\{e\}$ are left cells and that $E_{1}\cup E_{2}$, $F$ and $\{e\}$ are two-sided cells (see Example \ref{Fs}). Thus we see that
$$W-A\cup B\cup C\cup E\cup F\cup \{e\}=D=D_{1}\cup D_{2}\cup D_{3}$$
is a union of left and two-sided cells. For $1\leq i\leq 3$ we have $D\cap \cR^{\{s_{i}\}}=D_{i}$ thus $D_{i}$ is a union of left cells. Since $D_{i}$ is included in a left cell it follows that $D_{i}$ is a left cell.
Using Lemma \ref{itsc}, one can easily check that
$$D_{1}\cup D_{2}\cup D_{3}$$
is a two-sided cell.

\subsection{Left and two-sided order}
\label{leftorder}
\begin{Th}
The partial order induced by $\leq_{L}$ on the left cells can be described by the following Hasse diagram
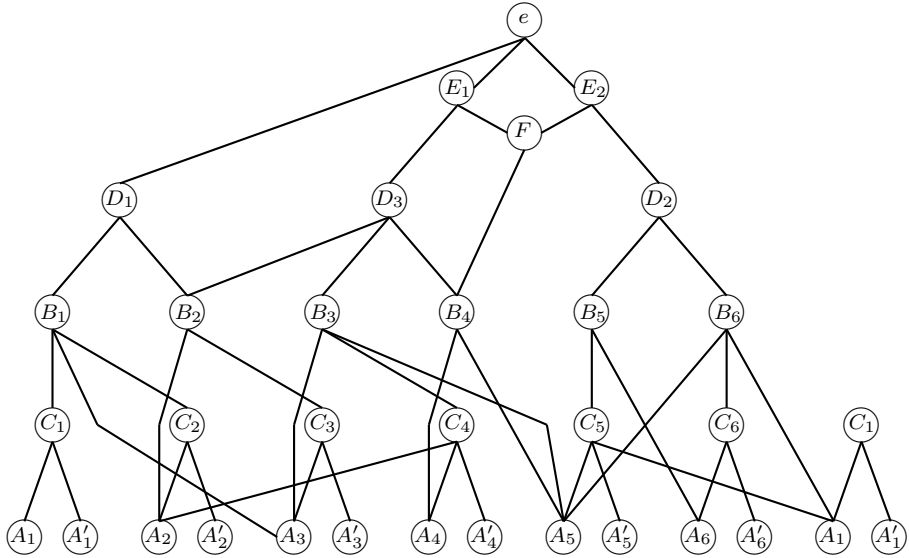
\begin{center}
\begin{pspicture}(0,0)(13,8)
\psset{unit=0.7cm}
\psset{unit=0.85pt}
\put( 17, 20){\circle{15}}
\put( 10.5, 18){\scriptsize{$A_{1}$}}
\put( 42, 20){\circle{15}}
\put( 35.5, 18){\scriptsize{$A'_{1}$}}

\put( 77, 20){\circle{15}}
\put( 70.5, 18){\scriptsize{$A_{2}$}}
\put( 102, 20){\circle{15}}
\put( 95.5, 18){\scriptsize{$A'_{2}$}}

\put( 137, 20){\circle{15}}
\put( 130.5, 18){\scriptsize{$A_{3}$}}
\put( 162, 20){\circle{15}}
\put( 155.5, 18){\scriptsize{$A'_{3}$}}

\put( 197, 20){\circle{15}}
\put( 190.5, 18){\scriptsize{$A_{4}$}}
\put( 222, 20){\circle{15}}
\put( 215.5, 18){\scriptsize{$A'_{4}$}}

\put( 257, 20){\circle{15}}
\put( 250.5, 18){\scriptsize{$A_{5}$}}
\put( 282, 20){\circle{15}}
\put( 275.5, 18){\scriptsize{$A'_{5}$}}

\put( 317, 20){\circle{15}}
\put( 310.5, 18){\scriptsize{$A_{6}$}}
\put( 342, 20){\circle{15}}
\put( 335.5, 18){\scriptsize{$A'_{6}$}}

\put( 29.5, 70){\circle{15}}
\put( 24, 68){\scriptsize{$C_{1}$}}

\put( 89.5, 70){\circle{15}}
\put( 84, 68){\scriptsize{$C_{2}$}}

\put( 149.5, 70){\circle{15}}
\put( 144, 68){\scriptsize{$C_{3}$}}

\put( 209.5, 70){\circle{15}}
\put( 204, 68){\scriptsize{$C_{4}$}}

\put( 269.5, 70){\circle{15}}
\put( 264, 68){\scriptsize{$C_{5}$}}

\put( 329.5, 70){\circle{15}}
\put( 324, 68){\scriptsize{$C_{6}$}}

\put( 29.5, 120){\circle{15}}
\put( 24, 117){\scriptsize{$B_{1}$}}

\put( 89.5, 120){\circle{15}}
\put( 84, 117){\scriptsize{$B_{2}$}}

\put( 149.5, 120){\circle{15}}
\put( 144, 117){\scriptsize{$B_{3}$}}

\put( 209.5, 120){\circle{15}}
\put( 204, 117){\scriptsize{$B_{4}$}}

\put( 269.5, 120){\circle{15}}
\put( 264, 117){\scriptsize{$B_{5}$}}

\put( 329.5, 120){\circle{15}}
\put( 324, 117){\scriptsize{$B_{6}$}}

\put( 59.5, 170){\circle{15}}
\put( 53, 168){\scriptsize{$D_{1}$}}

\put( 179.5, 170){\circle{15}}
\put( 173, 168){\scriptsize{$D_{3}$}}

\put( 299.5, 170){\circle{15}}
\put( 293, 168){\scriptsize{$D_{2}$}}

\put( 239.5,200){\circle{15}}
\put( 235, 198){\scriptsize{$F$}}

\put( 209.5, 220){\circle{15}}
\put( 204, 217){\scriptsize{$E_{1}$}}

\put( 269.5, 220){\circle{15}}
\put( 264, 217){\scriptsize{$E_{2}$}}

\put( 239.5,250){\circle{15}}
\put( 237, 249){\scriptsize{$e$}}

\put( 377, 20){\circle{15}}
\put( 370.5, 18){\scriptsize{$A_{1}$}}
\put( 402, 20){\circle{15}}
\put( 395.5, 18){\scriptsize{$A'_{1}$}}

\put( 389.5, 70){\circle{15}}
\put( 384, 68){\scriptsize{$C_{1}$}}

\psline(377,27,5)(389.5,62.5)
\psline(402.5,27,5)(389.5,62.5)

\psline(17,27,5)(29.5,62.5)
\psline(77,27,5)(89.5,62.5)
\psline(137,27,5)(149.5,62.5)
\psline(197,27,5)(209.5,62.5)
\psline(257,27,5)(269.5,62.5)
\psline(317,27,5)(329.5,62.5)

\psline(42,27,5)(29.5,62.5)
\psline(102,27,5)(89.5,62.5)
\psline(162,27,5)(149.5,62.5)
\psline(222,27,5)(209.5,62.5)
\psline(282,27,5)(269.5,62.5)
\psline(342,27,5)(329.5,62.5)

\psline(29.5,112.5)(29.5,77.5)

\psline(29.5,112.5)(89.5,77.5)

\psline(29.5,112.5)(49.5,70)
\psline(49.5,70)( 129.5, 20)

\psline(269.5,112.5)(269.5,77.5)
\psline(329.5,112.5)(329.5,77.5)

\psline(89.5,112.5)(149.5,77.5)

\psline(89.5,112.5)(77,70)
\psline(77,70)(77,27,5)

\psline(149.5,112.5)(137,70)
\psline(137,70)(137,27,5)

\psline(209.5,112.5)(197,70)
\psline(197,70)(197,27,5)

\psline( 209.5, 112.5)(257,27,5)

\psline(269.5,112.5)(317,27,5)


\psline(329.5,112.5)(377,27.5)

\psline(329.5,112.5)(257,27,5)


\psline(209.5,62.5)(77,27,5)

\psline(269.5,62.5)(377,27,5)

\psline(149.5,112.5)(209.5,77.5)

\psline(149.5,112.5)(249.5,70)
\psline(249.5,70)(257,27,5)

\psline(59.5, 162.5)(29.5,127.5)

\psline(59.5, 162.5)(89.5,127.5)

\psline(179.5, 162.5)(149.5,127.5)

\psline(179.5, 162.5)( 209.5, 127.5)


\psline( 299.5, 162.5)(269.5,127.5)

\psline( 299.5, 162.5)(329.5,127.5)


\psline( 239.5,192.5)( 209.5, 127.5)

\psline( 209.5, 212.5)( 232,200)

\psline( 209.5, 212,5)(179.5, 177.5)


\psline( 269.5, 212.5)( 247,200)

\psline( 269.5, 212.5)( 299.5, 177.5)


\psline(239.5,242,5)( 217, 220)

\psline(239.5,242.5)( 262, 220)

\psline(239.5,242,5)(59.5, 177.5)

\psline(179.5, 162.5)(89.5,127.5)
\end{pspicture}
\end{center}
\end{Th}
\smallskip
\begin{proof}
Most of the relations can be deduced using the fact that for $s\in S$ and $w\in W$, if $sw>w$ then $sw\leq_{L} w$. For instance, for all $1\leq i\leq 6$ we have $A_{i}\leq_{L} C_{i}$ and $A'_{i}\leq_{L} C_{i}$.\\
Some of the relations require some explicit computations, we refer to \cite{webpage} for details. The fact that there is no other links comes from the last two sections, where we have determined many left ideals of $W$. Recall that in \cite{jeju2}, it is shown that $A_{i}$ and $A'_{i}$ are left ideals of $W$.

\end{proof}

\begin{Th}
Let $T=D$ or $T=F=\{s_{2}s_{3}s_{2}\}$. Then the partial order induced by $\leq_{LR}$ on the two-sided cells is as follows
$$A\leq C\leq B\leq T\leq E\leq \{e\}$$
and $D$ and $F$ are not comparable.\\
\end{Th}
\begin{proof}
This is easily checked.
\end{proof}
Using the explicit decomposition of $\tilde{G}_{2}$ in our case, we can check some of Lusztig's conjectures (see \cite[Chap. 14]{bible}). For instance
$${\bf P14}. \quad \text{For any $z\in W$, we have $z\sim_{LR} z^{-1}$}$$
is certainly true. The following statement can be easily deduced from {\bf P4} and {\bf P9}
$$x\leq_{L} y \quad \text{and} \quad x\sim_{LR} y \quad \Longrightarrow\quad x\sim_{L} y.$$
This can be easily checked from the partial left order on the left cells. Indeed, there is no relations between two left cells lying in the same two-sided cell.\\

\noindent{\bf Acknowledgment}.
I would like to thank C\'edric Bonnaf\'e for pointing out a gap in the proof of Theorem \ref{tfix}. I would also like to thank the referee for some very useful comments.

\bibliography{biblio1}

\end{document}